\documentclass[12pt]{article}
\usepackage[utf8]{inputenc}
\usepackage[english]{babel}
\usepackage{graphicx}
\usepackage{amssymb}
\usepackage{amsthm}

\usepackage{subcaption}
\usepackage[font=footnotesize]{caption}

\usepackage{tikz}
\usetikzlibrary{matrix}
\usepackage[left=1.5in,right=1in]{geometry}

\usepackage[T1]{fontenc}

\usepackage[shadow, colorinlistoftodos,textsize=tiny]{todonotes}
\setuptodonotes{fancyline, backgroundcolor=gray!10,bordercolor=gray}
\newcommand{\missingref}[1][?]{\todo[color=gray!30]{Reference (#1)}}
\usepackage{authblk}
\usepackage{hyperref}
\hypersetup{
	hidelinks,
  colorlinks  = true,    
  urlcolor    = blue,    
  linkcolor   = blue,    
  citecolor   = blue,      
  pdfauthor   = {Isabel Hubard, Elías Mochán, Antonio Montero },%
  pdfsubject  = {Maniplexes},%
  pdftitle    = {Voltage operations on maniplexes},  %
}

\usepackage[capitalise,noabbrev, nameinlink, poorman]{cleveref}

 \newcommand{\elias}[1]{\todo[color=red!25]{\textbf{Elías:} #1}}

\newtheorem{teo}{Theorem}[section]
\newtheorem*{Mteo}{Main Theorem}
\newtheorem{coro}[teo]{Corollary}
\newtheorem{prop}[teo]{Proposition}

\newtheorem{lema}[teo]{Lemma}
\theoremstyle{definition}
\newtheorem{Def}[teo]{Definition}

\newtheorem{preg}{Question}

\newenvironment{dem}[1][Proof]{\begin{trivlist}
  \item[\hskip \labelsep {\bfseries #1}]}{\end{trivlist}}

\theoremstyle{remark}
\newtheorem{obs}[teo]{Remark}

\newcommand{\F}{\mathcal{F}}

\newcommand{\md}{{\rm \ mod\ }}
\newcommand{\floor}[1]{\left \lfloor #1 \right \rfloor}

\newcommand{\p}{\mathcal{P}}
\newcommand{\Q}{\mathcal{Q}}
\newcommand{\G}{\mathcal{G}}
\newcommand{\M}{\mathcal{M}}
\newcommand{\T}{\mathcal{T}}

\newcommand{\U}{\mathcal{U}}

\newcommand{\ol}[1]{\overline{#1}}
\newcommand{\os}[1]{\widetilde{#1}}

\newcommand{\gen}[1]{\langle #1 \rangle} 
\newcommand{\Aut}{\Gamma}

\newif{\ifEjemploPrimero}
\newif{\ifDetallesGeneradores}
\newif{\ifCondsOrdenRaro}
\newif{\ifDetallesEjIntProp}
\newif{\ifEjemploConstruccion}
\newif{\ifEjemploPropsTresOrb}
\newif{\ifEjemploManiplex}

\title{All polytopes are coset geometries: characterizing automorphism groups of $k$-orbit abstract polytopes
}
\author{Isabel Hubard}
\affil{Institute of Mathematics, Universidad Nacional Aut\'onoma de
  M\'exico (IM UNAM), 04510 Mexico City, Mexico\thanks{\tt email:isahubard@im.unam.mx}}
\author{Elías Mochán}
\affil{Department of Mathematics, Northeastern University, 02115
  Boston, USA\thanks{\tt email:j.mochanquesnel@northeastern.edu}}
\begin{document}
\maketitle

\abstract{Abstract polytopes generalize the classical notion of convex polytopes to more general combinatorial structures. 
The most studied ones are regular and chiral polytopes, as it is well-know they can be constructed as coset geometries from their automorphism groups. This is also known to be true for 2- and 3- orbit 3-polytopes.
In this paper we show that every abstract n-polytope can be constructed as a coset geometry.
This construction is done by giving a characterization, in terms of generators, relations and intersection conditions, of the automorphism group of a $k$-orbit polytope with given symmetry type graph.
Furthermore, we use these results to show that for all $k\neq 2$, there exist $k$-orbit $n$-polytopes with Boolean groups (elementary abelian $2$-groups) as automorphism group, for all $n\geq 3$.}

\section{Introduction}
In the 1970s several ideas extend the geometric study of convex polytopes to  generalize them from different points of view:
while Tits studied incidence systems, Coxeter focused on tessellations of manifolds and Gr\"unbaum proposed to move away from spherical tiles.
In the early 1980s Danzer and Schulte put several of these ideas together to start the study of {\em incidence polytopes}, now called {\em abstract polytopes}.
An abstract polytope is a ranked, partially ordered set which generalizes the face lattice of convex polytopes and tessellations. Thus, abstract polytopes generalize the classical notion of convex polytopes and tessellations to more general combinatorial structures.

The degree of symmetry of an abstract polytope is measured by counting the number of flag
orbits under the action of its automorphism group, where a flag is a maximal chain of the partial order.
Abstract polytopes can then be classified in terms of their so-called symmetry type graph (\cite{SymType}), which encodes all the information of the local configuration of flags with respect to the automorphism group.

Polytopes with only one flag-orbit
are called {\em regular} and  are  the most studied ones.
Regular polytopes have maximal degree of symmetry and in particular their groups are generated by involutions, often called ``(abstract) reflections''.
The
book~\cite{ARP} is the standard reference and it is devoted exclusively
to the study of abstract regular polytopes.
There are $2^{n}-1$ classes of n-polytopes with 2 flag orbits, each of them corresponding to one symmetry type graph.
Among them is the class of {\em chiral} polytopes (\cite{QuiralesEgonAsia}), which have no ``reflectional'' symmetry, but have maximal ``rotational'' one.

The study of coset geometries goes back to Tits (\cite{Tits61}).
The ideas behind this concept are to construct incidence structures using groups, and in particular were developed by Tits in connection to Coxeter groups.
It is well-known that regular and chiral polytopes, as well as two-orbit polyhedra can be seen as coset geometries (see \cite{EgonPhd}, \cite{QuiralesEgonAsia} and \cite{2OrbPol1}); this characterization of their automorphism groups constitute the most important tool to study them.
In a recent paper (\cite{3Orb}), we showed that also 3-orbit polyhedra can be seen as cosets and used this to construct 3-orbit polyhedra from symmetric groups.
%

The purpose of this paper is to show the following theorem.
\begin{Mteo}
Every abstract polytope can be constructed as a coset geometry.
\end{Mteo}
 To show this, we characterize the automorphism group of an abstract polytope with a given symmetry type graph, in terms of generators and relations, as well as some intersection conditions on some subgroups and cosets.

Throughout the paper, we work with the flag graph of a polytope, as opposed to  the partial order. In fact, we shall often work with {\em maniplexes}, that is, colored graphs that generalize flag graphs of polytopes (\cite{Maniplexes}).

As pointed out before, symmetry type graphs are a great tool to classify polytopes (and maniplexes) in terms of their automorphism groups.
We will rely on them heavily in our study.
There are some (necessary)
conditions a graph must satisfy to be the symmetry type graph of a
maniplex or polytope.
Graphs satisfying such conditions are called
{\em premaniplexes} (or admissible graphs) and will be properly defined in Section~\ref{s:STG}.
We believe our main theorem will shed light to the study of $k$-orbit polytopes, at the time that opens the gate to study one of the main problems in the area (listed as Problem 12 in~\cite{k-orb}):

\begin{preg}\label{q:premaniplexasSTG}
  Given a premaniplex, does it exist a polytope (or maniplex) having such premaniplex as its symmetry type graph.
\end{preg}

The problem of finding polytopes or maniplexes with a given symmetry type graph is
in general very difficult, as one can note, for example, by looking
at the history of chiral polytopes:
although it was back in 1991~\cite{QuiralesEgonAsia}, when Schulte and Weiss studied chiral
polytopes and classified their automorphism groups in terms of
generators and relations, it took almost $20$ years to have a construction showing that such polytopes existed for all ranks $n>3$ (see \cite{Quirales}).


Very little is know about other particular instances of Question~\ref{q:premaniplexasSTG}.
It is known~\cite{SymType} that every
premaniplex with 3 vertices with rank $n\geq 3$ is the symmetry
type graph of a polytope,
and in~\cite{2OrbMani} Pellicer, Potočnik and Toledo construct 2-orbit maniplexes, but it is not known if they are polytopal
(i.e. the flag graph of a polytope).

This paper is organized as follows. In Section~\ref{sec:Polytopes-and-maniplexes} we give the basic concepts from the theory of abstract polytopes as well as maniplexes, and formally introduce the concepts of symmetry type graphs and premaniplexes. We also state a relaxed version of Question~\ref{q:premaniplexasSTG}, that we answer later, in Section~\ref{s:IntProp}.
We start Section~\ref{sec:voltage} by going over the main tool used in this paper: voltage graphs. Then,
by using them, in Section~\ref{s:VoltMani} we construct a maniplex $\M$ from
a premaniplex $X$ and a group $G$ (satisfying some
conditions).
Then $G$ will act on $\M$ by automorphisms and the
quotient of $\M$ by the action of $G$ will be $X$.
This means that $\M$ will have symmetry type graph $X$ if and only if every automorphism of $\M$ is represented by the action of an element of $G$.

In Section~\ref{s:IntProp} we use the construction of Section~\ref{sec:voltage} as well as the results from \cite{PolyMani} to characterize in terms of generators and relations the groups that are automorphisms groups of a polytope with a given symmetry type graph. We do so by showing that automorphism groups of polytopes must satisfy certain ``intersection conditions'' for some subgroups and cosets, that depend on the symmetry type graph.
In other words, we give an algebraic test for the group $G$ of Section~\ref{s:VoltMani}  that tells us if
the constructed maniplex $\M$ is polytopal or not, thus translating
the problem of finding polytopes with a given symmetry type graph to a
group-theoretic one.
This gives an answer to Problem 1 of \cite{k-orb}.
In Section~\ref{s:Construction} we give the proof to the Main Theorem by constructing a polytope as a coset geometry from a voltage group.
This gives an answer to Problem 2 of \cite{k-orb}.

We finish the paper in Section~\ref{s:Caterpillars} by using the above result to construct (degenerated) $k$-orbit $n$-polytopes with with Boolean groups (elementary abelian $2$-groups) as automorphism group, for all $k, n\geq 3$.
For this, we define {\em caterpillars} as premaniplexes having exactly one generating three and studying their coverings in order to avoid the possible extra symmetry that might happen when one uses our construction to obtain polytopes from groups.

\section{Abstract polytopes and maniplexes}\label{sec:Polytopes-and-maniplexes}

In this section we shall give the basic definitions and properties of abstract polytopes and maniplexes, and some relations between them. For more details, we refer the reader to \cite{PolyMani}, \cite{ARP} and \cite{Maniplexes}.

A partially ordered set is said to be {\em flagged} if it has
a (unique) least and a (unique) greatest
element and each maximal chain, called {\em a flag}, has the same
finite cardinality.
As all flags have the same cadrinality, say $n+2$, flagged posets naturally admit an order-preserving
function, the {\em rank},  from the poset to the set $\{-1,0,1,\ldots,n\}$.
The rank function allows us to talk about flag-adjacencies: given two flags $\Phi$ and $\Psi$ of a flagged poset, they are said to be $i$-adjacent if they satisfy to differ only in the element of rank $i$.

An {\em (abstract) $n$-polytope} (also called an {\em (abstract)
  polytope of rank
  $n$}) is a flagged poset $\p$ in which the following
conditions hold:
\begin{itemize}
\item {\em Diamond condition:} for each flag $\Phi$ of $\p$ and each $i \in \{0,1,\dots, n-1\}$, there exists a unique $i$-adjacent flag to $\Phi$. 
  \item {\em Strong flag connectedness:} for any two flags $\Phi$ and
    $\Psi$ of $\p$, 
    there exists a sequence of adjacent flags connecting $\Phi$ to $\Psi$,
    such that all the flags in the sequence contain the faces of the
    intersection $\Phi \cap \Psi$.
\end{itemize}


The elements of rank $i$ in a polytope are called {\em $i$-faces}.
Given a flag $\Phi$, we often denote its $i$-face as
$\Phi_i$, and by $\Phi^i$ its (unique!)
$i$-adjacent flag.
Recursively, if $w$ is a word on
$\{0,\ldots,n-1\}$ we denote by $\Phi^{wi}$ the flag $(\Phi^{w})^i$.
It is straightforward to see that $(\Phi^i)^i=\Phi$, and that $\Phi_j^i=\Phi_j$ if and only if $i\neq j$.

Given a flagged poset $\p$ one can define its {\em flag graph}
as the graph $\G(\p)$ whose vertices are the flags of $\p$ and two
flags $\Phi$ and $\Psi$ are connected by an edge of color
$i\in\{0,1,\ldots,n-1\}$ if and only if they are $i$-adjacent.
The flag graph of an $n$-polytope is an {\em $n$-maniplex}, that is,
an $n$-regular connected simple graph with a proper edge coloring with colors
$\{0,\ldots,n-1\}$ such that if $i$ and $j$ are two colors satisfying that
$|i-j|>1$, then the graph induced by edges
of colors $i$ and $j$ is a disjoint union of 4-cycles.
However, not every maniplex is the flag graph of a polytope.

While abstract polytopes generalize classical polytopes to combinatorial structures,
maniplexes (introduced by Steve Wilson in 2012 \cite{Maniplexes}) generalize flag graphs of a polytope as well as of the flag graphs of maps on surfaces (that is, 2-cellular
embeddings of connected graph on a surface).
In order to unify our notation of abstract polytopes and maniplexes, when dealing with maniplexes we shall call {\em flags} to its vertices and say that two of them are $i$-adjacent if they are the vertices of an edge of colour $i$.

It follows from the definition of a maniplex that each flag is incident to exactly
one edge of each color.
Then, for each $i\in\{0,1,\ldots,n-1\}$ one
can define $r_i$ as the permutation of the set of flags that maps
each flag its $i$-adjecent one.
In other words $\Phi r_i=\Psi$ if and only if $\Phi$ and $\Psi$
are $i$-adjacent.

The permutations $r_i$, with $i\in \{0,1,\ldots,n-1\}$, are all
involutions with no fixed points and, by connectivity, they generate a group of permutations on the flags which acts
transitively on them.
Furthermore, if $|i-j|>1$ then $r_ir_j$ is also an involution with no
fixed points; thus, $r_i$ and $r_j$
commute.


Since $\Phi r_i$ is the flag $i$-adjacent to $\Phi$, it is convenient to
denote it by $\Phi^i$ and to follow the same recursive notation as before: $\Phi^{wi}=(\Phi^w)^i$ where $w$ is a word on
$\{0,1,\ldots,n-1\}$.

The group $\gen{r_0, r_1\ldots,r_{n-1}}$ is called the {\em monodromy} or
{\em connection} group of the maniplex $\M$, it shall be denoted by
$Mon(\M)$ and we shall call each of its elements a {\em
  monodromy}.
  If $w$ is a word on the
alphabet $\{0,1,\ldots,n-1\}$ we identify $w$ with the monodromy
$x\mapsto x^w$, that is ,the word $a_1a_2\ldots a_k$ is identified
with the monodromy $r_{a_1}r_{a_2}\ldots r_{a_k}$.

A {\em maniplex homomorphism} is a graph homomorphism that preserves the color of the edges.
Using the connectedness of maniplexes one can see that every maniplex homomorphism is determined by the image of one flag and that they are all surjective.
The notions of {\em isomorphism} and {\em automorphism} follow naturally.

 As with polytopes, we denote the automorphism group of a maniplex $\M$ by $\Gamma(\M)$.
 By definition, if $\gamma \in \Gamma(\M)$, then $(\Phi r_i)\gamma = (\Phi^i)\gamma = (\Phi\gamma)^i = (\Phi\gamma)r_i$, for all $i\in \{0,1, \dots, n-1\}$, implying that $\omega\gamma = \gamma \omega$ for all $\omega \in Mon(\M)$.

Note further that, since the action of $Mon(\M)$ is transitive, the action of $\Gamma(\M)$ is free (or semi-regular).
Of course, this is true for both  polytopes and maniplexes.






Let $\M$ be an $n$-maniplex. If $I\subset \{0,1,\ldots,n-1\}$, we
define $\M_I$ as the subgraph of $\M$ induced by the edges of colors
in $I$. If $i\in\{0,1,\ldots,n-1\}$, we use the symbol $\ol{i}$ to
denote the set $\{0,1,\ldots,n-1\}\setminus\{i\}$, and more generally,
if $K\subset \{0,1,\ldots,n-1\}$, we denote its complement by
$\ol{K}$. In particular $\M_{\ol{i}}$ is the subgraph of $\M$ obtained
by removing the edges of color $i$. We will use this notation also for
any graph with a coloring of its edges even if its not a maniplex.

In \cite{PolyMani} Garza-Vargas and Hubard describe how to recover a polytope $\p$ from its flag graph:
%
%
the elements of rank $i$ of $\p$ are the connected
components of $\M_{\ol{i}}$; the order is given as follows: given connected components  $F$ and $G$ of $\M_{\ol{i}}$ and $\M_{\ol{j}}$, respectively, we set $F<G$ if and only if $i<j$ and $F\cap G \neq \emptyset$


In \cite{PolyMani} it is proved that if $\M$ is any maniplex then
$\p(\M)$ is in fact a poset (actually, it is a flagged poset).

It is not difficult to see (\cite{PolyMani}) that $\p$ and $\p'$ are isomorphic polytopes if and only if their flag graphs $\G(\p)$ and $\G(\p')$ are isomorphic. This fact implies the following theorem, which can be interpreted as saying that all the information of the
polytope $\p$ is encoded in its flag graph.


\begin{teo}\cite{PolyMani}\label{t:GrafBand}
  Let $\p$ be a polytope and let $\M=\G(\p)$ be its flag graph. Then
  $\p$ is isomorphic (as a poset) to $\p(\M)$ and $\Gamma(\p) =
  \Gamma(\M)$.
\end{teo}

Theorem 5.3 of \cite{PolyMani} gives a characterization of {\em polytopal} maniplexes, that is, those maniplexes
that are isomorphic to the flag graph of some polytope.
Such characterization is given in terms of some path intersection properties of the maniplexes.


\begin{Def}\label{d:PIP}
  Let $\M$ be an $n$-maniplex. We say that $\M$ satisfies the {\em
    strong path intersection property} (or SPIP) if for every two
  subsets $I,J \subset \{0,1,\ldots,n-1\}$ and for any two flags $\Phi$ and
  $\Psi$, if there is a path $W$ from $\Phi$ to $\Psi$ using only darts of
  colors in $I$ and also a path $W'$ from $\Phi$ to $\Psi$ using only darts
  of colors in $J$, then there also exists a path $W''$ from $\Phi$ to
  $\Psi$ that uses only darts of colors in $I \cap J$.

  We say that $\M$ satisfies the {\em weak path intersection property} (or WPIP)
  if for any two flags $\Phi$ and $\Psi$ and for all $k,m \in
  \{0,1,\ldots,n-1\}$, whenever there is a path $W$ from $\Phi$ to $\Psi$
  with only darts of color in $[0,m]:=\{0,1,\ldots,m\}$ and a path $W'$
  from $\Phi$ to $\Psi$ with only darts of colors in
  $[k,n-1]:=\{k,k+1,\ldots,n-1\}$, then there is also a path $W''$
  from $\Phi$ to $\Psi$ with only darts of colors in $[k,m]:=\{k,k+1,\ldots,m\}$.
\end{Def}

\begin{teo}\label{t:PolyMani}
  \cite{PolyMani} Let $\M$ be a maniplex. Then the following conditions are all
  equivalent:
  \begin{itemize}
    \item $\M$ is polytopal.
  \item $\M$ satisfies the SPIP.
  \item $\M$ satisfies the WPIP.
  \end{itemize}

  In any of these cases $\p(\M)$ is a polytope whose flag graph is
  isomorphic to $\M$.
\end{teo}
\subsection{Premaniplexes and symmetry type graphs}\label{s:STG}
A \emph{$k$-orbit maniplex} is one with exactly $k$ flag orbits under its automorphism group.
When studying $k$-orbit maniplexes with $k>1$ one finds that it is convenient  to classify them in terms of the local structure of the flags.
For this reason, in~\cite{SymType}, Cunningham et al. introduce the concept of \emph{symmetry type graph}.

Given a maniplex $\M$ and a subgroup $G$ of the automorphism group of $\M$, the \emph{symmetry type graph of $\M$ with respect to $G$}, denoted either by $\T(\M,G)$ or by $\M/G$, is constructed as follows:
The vertex set of $\T(\M,G)$ is the set of flag orbits of $\M$ under the group $G$, and if $\Phi$ and $\Psi$ are $i$-adjacent in $\M$ we draw an edge of color $i$ between their orbits.
If $\Phi$ and $\Phi^i$ are in the same orbit under $G$, we draw a semi-edge of color $i$ at the vertex corresponding to that orbit.
Recall that a semi-edge is different from a loop in that it consists of only one dart which is inverse to itself, rather than two darts starting at the same vertex; less formally, a semi-edge is incident to the vertex once, while a loop is incident to the vertex two times.

When we speak about \emph{the symmetry type graph of $\M$}, we mean it with respect to $\Gamma(\M)$ and we simply write $\T(\M)$ instead of $\T(\M,\Gamma(\M))$.

If $\p$ is a polytope, \emph{the symmetry type graph of $\p$ (with respect to $G$)}  is defined as the symmetry type graph (with respect to $G$) of its flag graph.

If $X$ is the symmetry type graph of an $n$-maniplex, then it is a connected graph in which every vertex is incident to exactly one edge of each color in $\{0,1,\ldots,n-1\}$ and it satisfies that if $|i-j|>1$, the paths of length 4 that alternate between the colors $i$ and $j$ are closed.
However, it might not be a maniplex as it is not necessarily simple.
We will call such a graph a \emph{$n$-premaniplex}.

\begin{figure}
    \centering
    \includegraphics[width=0.5\textwidth]{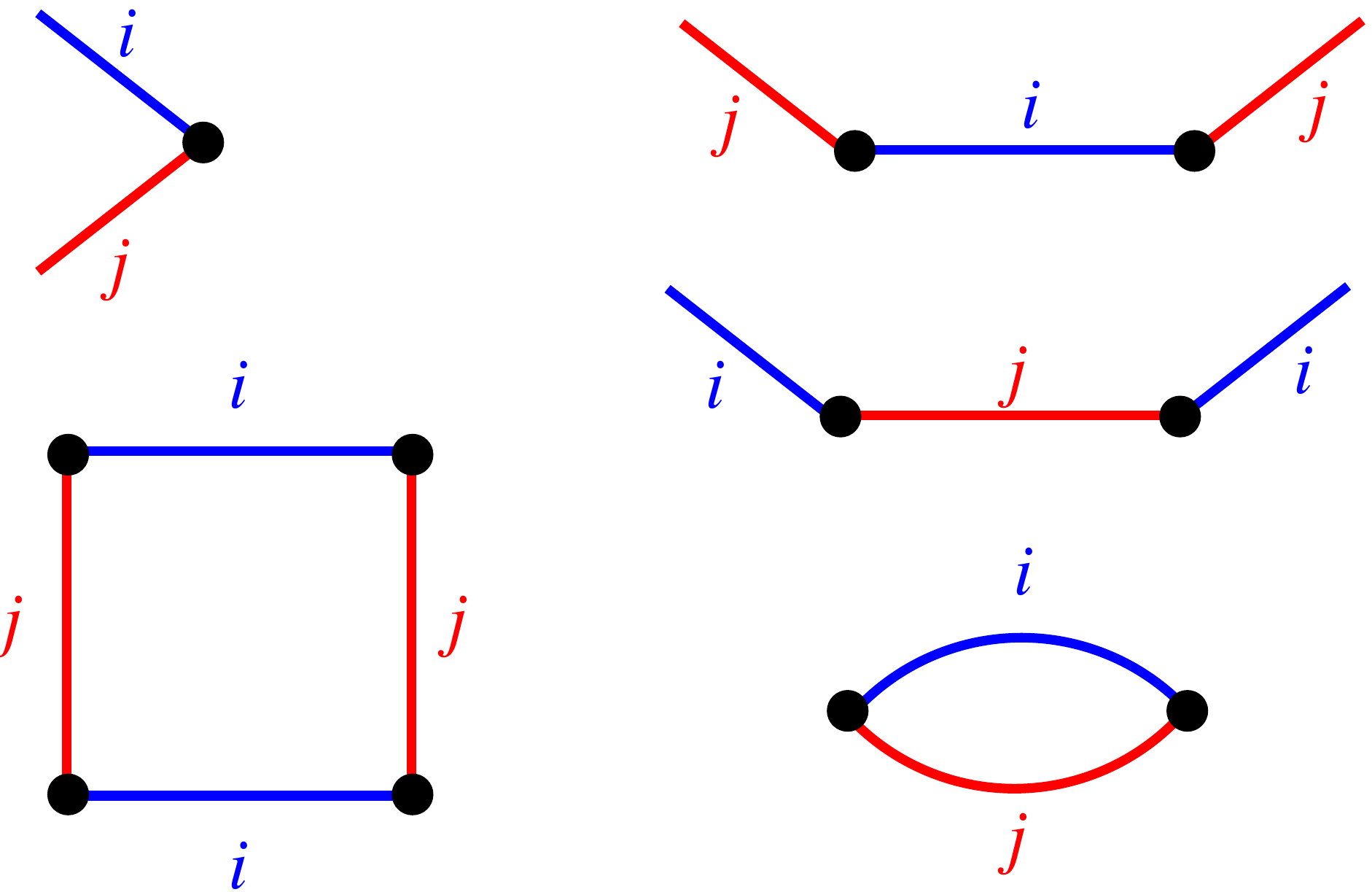}
    \caption{Components of two non-consecutive colors in a premaniplex.}
    \label{f:ijComponentes}
\end{figure}

If $X$ is a premaniplex and $i,j\in\{0,1,\ldots,n-1\}$ are non-consecutive, the connected components of the subgraph of $X$ induced by the edges of colors $i$ and $j$ are not necessarily 4-cycles. In fact they can be any quotient of a 4-cycle, as ilustrated in Figure~\ref{f:ijComponentes}.

The notions of \emph{homomorphism, isomorphism, automorphism} and \emph{monodromy group} from maniplexes can all be easily extended to premaniplexes as well.

The natural projection $p:\M\to\T(\M,G)$ is, of course, a homomorphism (of premaniplexes).

When studying polytopes (or maniplexes) and their symmetry type graphs two natural questions occur:
\begin{preg}\label{q:SymType}
    Given a premaniplex $X$, is there a polytope (or maniplex) whose symmetry type graph is $X$?
\end{preg}

\begin{preg}\label{q:IntProp}
    Given a premaniplex $X$, what conditions must a group $G$ satisfy so that there is a polytope (or maniplex) $\p$ such that $\T(\p,G)\cong X$?
\end{preg}


In \cref{s:IntProp} we give a complete answer to \cref{q:IntProp}.
\cref{q:SymType} remains as a hard question, as even if we find a polytope $\p$ and a group $G$ such that $\T(\p,G)\cong X$, it may still happen that $G$ is a proper subgroup of $\Gamma(\p)$.
However, in \cref{s:Caterpillars} we give an infinite family of premaniplexes that are in fact symmetry type graphs of polytopes.
\section{Voltage graphs}\label{sec:voltage}

The projection $p:\M\to\T(\M,G)$ is an example of what is called a \emph{regular covering projection} in graph theory (see~\cite{VoltsOrig} or~\cite{Voltajes} for more details).
Given a regular covering projection $p:\os{X}\to X$, one may recover the graph $\os{X}$ from the graph $X$ using what is known as a \emph{voltage assignment}.

A \emph{path} in a graph is a finite sequence of darts $W=d_1d_2\ldots,d_k$ such that the dart $d_{i+1}$ starts at the endpoint of the dart $d_i$.
The \emph{startpoint of $W$} is the starting point of $d_1$ (denoted $I(d_1)$), and the \emph{endpoint of $W$} is the endpoint of $d_k$ (denoted $T(d_k)$).
If $x$ and $y$ are the startpoint  and endpoint of $W$, respectively, we say that $W$ \emph{goes from $x$ to $y$} and we write this as $W:x\to y$.
If the endpoint and startpoint of a path $W$ are the same vertex $x$ we say that $W$ is a \emph{closed path based at $x$}.
We also consider that for every vertex $x$ there is an empty closed path 
based at $x$.

If a path $W$ ends at the startpoint of a path $V$, we may define the product $WV$ as their concatenation.

Two paths $W$ and $W'$ with the same startpoint and endpoint are said to be \emph{homotopic} if one can transform $W$ into $W'$ by a finite sequence of the following operations:
\begin{itemize}
    \item Inserting two consecutive inverse darts at any point, that is $$d_1d_2\ldots d_i d_{i+1}\ldots d_k\mapsto d_1\ldots d_i d d^{-1} d_{i+1}\ldots d_k,$$ where $I(d)=T(d_i)$;
    \item Deleting two consecutive inverse darts at any point, that is $$d_1\ldots d_i d d^{-1} d_{i+1}\ldots d_k \mapsto d_1d_2\ldots d_i d_{i+1}\ldots d_k;$$
\end{itemize}
In this case we write $W\sim W'$.

 It is easy to see that homotopy is an equivalence relation and that if $W\sim W'$ and $V\sim V'$, then $WV\sim W'V'$.
 Therefore, we can think of the product of two homotopy classes of paths.
 The set of all homotopy classes of paths in a graph $X$ with this operation is called \emph{the fundamental groupoid of $X$} and it is denoted by $\Pi(X)$.
 We will often speak of a "path $W$ in $\Pi(X)$", but the reader should keep in mind that we are actually referring to its homotopy class.

 The subset of $\Pi(X)$ consisting of all the (homotopy classes of) closed paths based at a vertex $x$ forms a group known as \emph{the fundamental group of $X$ (based at $x$)} and it is denoted by $\Pi^x(X)$.

 A \emph{voltage graph} is a pair $(X,\xi)$ where $X$ is a graph and $\xi$ is a groupoid antimorphism from $\Pi(X)$ to a group $G$. In this case we say that $\xi$ is a \emph{voltage assignment (on $X$)} and we call $G$ the \emph{voltage group (of $\xi$)}.
 The element $\xi(W)$ is called \emph{the voltage of $W$}.

 Note that a voltage assignment is completely determined by the voltages of the darts of the graph, as the voltage $\xi(W)$ of a path $W=d_1d_2\ldots d_{k-1}d_k$ is simply $\xi(d_k)\xi(d_{k-1})\ldots \xi(d_1)$.

 Given a voltage graph $(X,\xi)$ with voltage group $G$, we can construct the \emph{derived graph} $X^\xi$ as follows:
 \begin{itemize}
     \item The vertex set is $V\times G$ where $V$ is the vertex set of $X$.
     \item The dart set is is $D\times G$ where $D$ is the dart set of $X$.
     \item The dart $(d,g)$ starts at the vertex $(x,g)$ and ends at the vertex $(x,\xi(d)g)$.
 \end{itemize}
 This is an undirected graph, as the inverse of the dart $(d,g)$ is the dart $(d^{-1},\xi(d)g)$.

 Given a path $W$  starting at a vertex $x$ and given an element $g$ in the voltage group, there is a unique path $\os{W}$ in $X^\xi$ that starts at $(x,g)$ that projects to $W$.
 The path $\os{W}$ is called a \emph{lift of $W$} and it is easy to see that it ends at $(x,\xi(W)g)$ (for details see~\cite{Voltajes}).

 In our case, we will work with graphs that have a coloring of its edges (and therefore, darts), so we will define that the color of the dart $(d,g)$ is the same as the color of $d$.

 It is known (see~\cite{VoltsOrig} or ~\cite{Voltajes}) that if $p:\os{X}\to X$ is a regular covering, there is a voltage assignment $\xi$ such that $\os{X}$ is isomorphic to $X^\xi$.

\subsection{Voltage graphs that give maniplexes as derived graphs}\label{s:VoltMani}

Of course, we want to use voltage graphs to obtain maniplexes (and polytopes). In this section we shall give necessary and sufficient conditions on a voltage assignment $\xi$ on a premaniplex $X$ so that $X^\xi$ is a maniplex.

In \cite[Section 5]{SymType} the authors find a set of distinguished generators for the automorphism group of a maniplex $\M$ with a given symmetry type graph $X$.
It is easy to see that such distinguished generators can be thought of as voltages assigned to $X$ (or more precisely to $\Pi(X)$) so that the derived graph is $\M$.
To do this, one must choose a spanning tree $T$ on $X$ and assign trivial voltage to all its darts. 
If we want to start with a premaniplex $X$, we need to be careful with the way we assign voltages so that the derived graph is indeed a maniplex.

\if
If the STG of $\M$ is $X$, there is a method to assign voltages to the remaining darts in such a way that the derived graph is isomorphic to $\M$.
It is known (\cite{SymType}\missingref{agregar resultado exacto}) that
the voltages of these darts have to generate the whole voltage group for the derived graph to be connected. So we call the non-trivial voltages of the darts
\fi

Let $X$ be a premaniplex with fundamental groupoid $\Pi(X)$, and let $\xi:\Pi(X)\to \Gamma$ be a voltage assignment, for some group $\Gamma$.
We want to find the conditions on $\Gamma$ and $\xi$ that ensure that $X^\xi$ is actually a maniplex,
but before doing so, let us assume that there is a spanning tree $T$ of $X$ with trivial voltage in all its darts.

First we want $X^\xi$ to be connected.
It is known (see~\cite{Voltajes}) that in order for $X^\xi$ to be connected, $\xi(D)$ must generate $\Gamma$, where $D$ denotes the set of darts of $D$.

Next, $X^\xi$ must be a simple graph.
Thus, it must not have semi-edges nor multiple edges.
Note that a semi-edge of $X^\xi$ that starts at a vertex $(x,\gamma)$ ends in
$(x,\xi(e)\gamma) = (x,\gamma)$, where $e$ is a semi-edge of $X$.
This implies that if the voltage of every semi-edge of $X$ is not trivial, we avoid semi-edges in $X^\xi$.
Since $\xi$ is an antimorphism we should have that $\xi(e)=\xi(e^{-1}) = \xi(e)^{-1}$, implying that the voltage of a semi-edge must have order two.

%
To avoid multiple edges, we need to avoid different darts with the same initial and terminal vertices;
suppose $X^\xi$ has two parallel darts
$(d,\sigma)$ and $(d',\tau)$.
Since both darts start at the same
vertex, we have that $(I(d),\sigma)=(I(d'),\tau)$, so $I(d)=I(d')$ and $\sigma =
\tau$. The common end-point of $(d,\sigma)$ and $(d',\tau)$ could be
written as $(y,\xi(d)\sigma)$
or $(z,\xi(d')\sigma)$, where $y$ is the end-point of $d$ and $z$ the
end-point of $d'$. The fact that these two are the same means
that $y=z$ and $\xi(d) = \xi(d')$. So $(d,\sigma)$ and
$(d',\sigma)$ are parallel darts in $X^\xi$ if and only if $d$ and $d'$ are
parallel darts in $X$ with the same voltage. Thus, $X^\xi$ has no parallel
darts if and only if no pair of parallel darts in $X$ has equal voltages.

Finally, we want to ensure that if $|i-j|>1$, the paths of length 4 in
$X^\xi$ that alternate colors between $i$ and $j$ are closed. Let
$\os{W}$
be one of these paths.
Projecting $\os{W}$ to $X$ we get a path $W$ in
$X$ of length 4 that alternates colors between $i$ and $j$, and
since $X$ is a premaniplex we know that $W$ is closed.
Suppose $W$ starts at a vertex $x$.
Then $\os{W}$ goes from a vertex of the form
$(x,\gamma)$ to $(x,\xi(W)\gamma)$. So $\os{W}$ is closed if and only
if $\xi(W)\gamma = \gamma$, or in other words, $W$ has trivial
voltage.
Summarizing this discussion, we arrive to the following lemma:

\begin{lema}\label{l:VoltMani}
  Let $X$ be a premaniplex and let $\xi:\Pi(X) \to \Gamma$ be a
  voltage assignment with a spanning tree $T$ of trivial voltage on all
  its darts. Then $X^\xi$ is a maniplex if and only if
  \begin{enumerate}
  \item The set $\xi(D)$ generates $\Gamma$, where $D$ is the set of
    darts of $X$,
  \item  $\xi(d)$ has order exactly 2 when $d$ is a semi-edge,
  \item $\xi(d) \neq \xi(d')$ when $d$ and $d'$ are
    parallel darts, and
  \item if $|i-j|>1$ every (closed) path $W$ of length 4 that alternates between
    the colors $i$ and $j$ has trivial voltage.
  \end{enumerate}
\end{lema}

In Section~\ref{sec:int prop} we shall translate the conditions in Lemma~\ref{l:VoltMani} to relations and inequalities that the generators of a group $\Gamma$ need to satisfy to act on a maniplex with given symmetry type graph $X$.
Before doing so (in the next section),
Let us now take a closer look at the consequences of condition 4 of the above Lemma.

Condition 4 of Lemma~\ref{l:VoltMani} invites us to introduce a new concept of homotopy:
we say that two paths $W$ and $W'$  are \emph{maniplex-homotopic} if we can transform one into the other by a finite sequence of inserting or deleting pairs of inverse darts, as well as switching the colors of two consecutive darts with non-consecutive colors, that is $$d_1d_2\ldots d_i d_{i+1}\ldots d_k\mapsto d_1d_2\ldots d'_i d'_{i+1}\ldots d_k,$$ where $|c(d_i)-c(d_{i+1})|>1,$ $c(d'_i)=c(d_{i+1})$, $c(d'_{i+1})=d(d_i)$, $I(d'_i)=I(d_i)$ and $I(d'_{i+1})=T(d'_i)$.
Hence, a voltage assignment $\xi$ is well defined when applied to the maniplex-homotopy class of paths if and only if it satisfies Condition 4 of Lemma~\ref{l:VoltMani}.
From this point on, whenever we speak about homotopy, homotopy class, fundamental groupoid, etc. we will be thinking in terms of maniplex-homotopy.

One could use the group $\Gamma:=\gen{S\mid R}$ as the voltage group where
$S$ has a generator $\alpha_e$ for each edge $e$ not in the spanning
tree of $X$ and $R$ has one element $\alpha_e^2$ per each semi-edge and one element $\alpha_{e_4}\alpha_{e_3}\alpha_{e_2}\alpha_{e_1}$ for each path of length four alternating between two non-consecutive colors.
In fact every voltage
group that gives a maniplex should be a quotient of this group, or in
other words $\Gamma$ is the ``most general'' group we can use as a
voltage group to get a maniplex as the derived graph.
We know (see~\cite{Hartley_AllPtpsAreQuotients}) that $X=\T(\U,G)$ where $\U$ is the
universal polytope of rank $n$ and $G$ is some group. This means that
there is some voltage assignment $\xi$ on $X$ with voltage group $G$
such that $X^\xi$ is isomorphic to the flag graph of $\U$. Because of
the universality of $\U$ we get that $G$ and $\Gamma$ in fact are the
same.
In other words, if we use the most general group as our voltage
group we will always get the flag graph of the universal polytope as
the derived graph.

\section{Intersection properties and coset geometries}\label{sec:int prop}



In order to prove the Main Theorem, we need to characterize, in terms of generators and relations, the
groups $\Gamma$
that act by automorphisms on a polytope $\p$ in such a way that the symmetry type graph $\T(\G(\p),\Gamma)$ is isomorphic to a given premaniplex $X$.
We shall do this in order to be able to recover the polytope $\p$ as a coset geometry using the group $\Gamma$.

We start with the premaniplex $X$ and provide it with a 
voltage assignment $\xi$.
Recall that if $(X,\xi)$ is a voltage graph with voltage group $\Gamma$, then $X^\xi/\Gamma$ is isomorphic to $X$; conversely, if $\M/\Gamma$ is isomorphic to $X$ then there is a voltage assignment $\xi$ on $X$ such that $X^\xi$ is isomorphic to $\M$.
Hence, by characterizing the voltage assignments $\xi$  that satisfy that the derived graph $X^\xi$ is the flag graph of a polytope, we  determine the conditions that $\Gamma$ must satisfy to be the automorphism group of a polytope with symmetry type $X$.

\subsection{Voltage graphs and the path intersection property}\label{s:IntProp}

We have figured out how to construct a maniplex from a premaniplex via voltage assignments. It is then natural to ask: when is the obtained maniplex the flag graph of a polytope?
The answer, as we shall see in this section, is closely related to Theorem~\ref{t:PolyMani}. In fact, translating Theorem~\ref{t:PolyMani} to the setting of voltage assignments will give us
conditions that take the form of {\em intersection properties} that certain distinguished subgroups and some left cosets must satisfy. 
Given two vertices $x,y$ in $X$ and a set of colors $I\in\{0,1,\ldots,n-1\}$, let us denote by $\Pi^{x,y}_I(X)$ the set of (homotopy classes of) paths from $x$ to $y$ in $X$ that only use darts with colors in the set $I$.
So $\xi(\Pi^{x,y}_I(X)) $ denotes the set of voltages of all the paths of $X$ from $x$ to $y$ whose edges have colors in $I$.

\begin{teo}\label{t:IntProp-w}
  Let $X$ be a premaniplex and let $\xi:\Pi(X)\to \Gamma$ be a voltage assignment such that $X^\xi$ is a maniplex. Then $X^\xi$ is the flag graph of a polytope if and only if
  \begin{equation}\label{VolPIP}
    \xi(\Pi^{x,y}_I(X)) \cap \xi(\Pi^{x,y}_J(X)) = \xi(\Pi^{x,y}_{I \cap J}(X)),
  \end{equation}
  for all $I,J \subset \{0,\ldots,n-1\}$ and all vertices $x,y$ in $X$.
\end{teo}
\begin{dem}
  Start by assuming that $X^\xi$ is the flag graph of a polytope.
  Let $x$ and $y$ be vertices of $X$ and let
  $I,J\subset\{0,1,\ldots,n-1\}$.
  Consider two paths, $W\in \Pi^{x,y}_I(X)$ and $W'\in
  \Pi^{x,y}_J(X)$, with the same voltage, say $\alpha\in\Gamma$.
  When lifting $W$ and $W'$ in $X^\xi$, they lift to
  paths $\os{W}$ and $\os{W}'$, respectively, that go from $(x,1)$ to $(y,\alpha)$ (here 1 is the identity element of $\Gamma$) and satisfying that $\os{W}$ uses edges with colors in $I$ while $\os{W}'$ uses edges with colors in $J$.
  In fact, one would define $\os{W}$ and $\os{W}'$ as the paths that start at $(x,1)$ and follow the same sequence of colors as $W$ and $W'$ respectively.
  By Theorem~\ref{t:PolyMani} $X^\xi$ satisfies the
  SPIP, which implies that there is a path
  $\os{W}''$ from $(x,1)$ to $(y,\alpha)$ that uses only colors in $I \cap J$.
  Then, its projection $W'':= p(\os{W}'')$ is a path in $X$
  that goes from $x$ to $y$ that uses only colors in $I \cap J$ and
  has voltage $\alpha$.
  This proves that $\xi(\Pi^{x,y}_I(X)) \cap
  \xi(\Pi^{x,y}_J(X)) \subset \xi(\Pi^{x,y}_{I \cap J}(X))$. Since the other
  contention is given, equality~(\ref{VolPIP}) must hold.

  Now let us assume that equality~(\ref{VolPIP}) holds for all $I,J\subset
  \{0,1,\ldots,n-1\}$ and all vertices $x$ and $y$.
  Let $\os{W}$ and
  $\os{W}'$ be paths in
  $X^\xi$ from a vertex $(x,\gamma)$ to a vertex $(y,\tau)$. Let $I$ and $J$
  be the sets of colors of darts of $\os{W}$ and $\os{W}'$, respectively, and let $W:=p(\os{W})$
  and $W':=p(\os{W'})$ be the projections of the paths to $X$.
  Then, both $W$ and $W'$ go from $x$ to $y$, and $W$ is a path with colors in $I$, while $W'$ is a path with colors in $J$, that is, $W\in \Pi^{x,y}_I(X)$ and $W'\in
  \Pi^{x,y}_J(X)$.
  Furthermore, since $\os{W}$ and $\os{W}'$ start at $(x,\gamma)$ and finish at $(y,\tau)$, they both have voltage
  $\alpha:=\tau\gamma^{-1}$.
  By hypothesis, there exists a path $W''\in \Pi^{x,y}_{I \cap J}(X)$ that also has  voltage $\alpha$.
  Then $W''$ has a unique lift $\os{W}''$ which is a
  path in $X^\xi$ from $(x,\gamma)$ to $(y,\tau)$ and it uses darts
  of colors in $I \cap J$. This  proves that $X^\xi$ satisfies the
  SPIP and therefore, by Theorem~\ref{t:PolyMani}, it is the flag graph of a polytope. \qed
\end{dem}

Note that when $x=y$ the set $\xi(\Pi^{x,y}_I(X)) := \xi(\Pi^x_I(X))$ is a group,
since it is the image of a group under a groupoid
antimorphism. Actually, we shall find a set of distinguished generators
for the group $\xi(\Pi^x_I(X))$ in a similar way as  the distinguished
generators of the automorphism group of a polytope are found in
\cite{SymType}.
Recall that $X_I$ is the subgraph of
$X$ induced by the edges with colors in $I$ and that $X_I(x)$ is the
connected component of $X_I$ containing the vertex $x$.
To find the distinguished generators of $\xi(\Pi^x_I(X))$, fix a spanning
tree $T^x_I$ for $X_I(x)$. 
For each dart $d$ in $X_I(x)$ but not in
$T^x_I$
we get a cycle $C_d$ of the form $WdV$ where $W$ is the unique path
contained in $T^x_I$ from $x$ to the initial vertex of $d$, and $V$ is
the unique path contained in $T^x_I$ from the terminal vertex of $d$
to $x$.
Then, the set $\{C_d\}$, where $d$ runs among
the darts in $X_I(x)$ not in $T^x_I$, is a generating set for
$\Pi^x_I(X)$.
This imples that $\{\xi(C_d)\}$, where $d$ runs among the darts of $X_I(x)$ not in $T^x_I$, is a set of generators for
$\xi(\Pi^x_I(X))$.

Since $\xi$ is a voltage assignment, we might consider only one dart $d$ for each edge. By denoting by $W_y$ the unique path
contained in $T^x_I$ from $x$ to $y$, we can see that $\Pi_I^{x,y}(X) = \Pi_I^x(X)W_y$ (that is, a path from $x$ to $y$ can be written as a closed path starting and finishing at $x$, concatenated with $W_y$), which implies that $\xi(\Pi^{x,y}_I(X)) =
\xi(W_y)\xi(\Pi^x_I(X))$.
Therefore, all the
intersection properties can be given in terms of left cosets of the groups $\xi(\Pi^x_I)$, whose generators we already know.

Theorem~\ref{t:IntProp-w} gives an intersection property for each pair
of vertices $(x,y)$ and each two sets of colors $I,J \subset
\{0,1,\ldots,n-1\}$. If we prove an intersection property for the pair
$(x,y)$, by taking the inverse on both sides we get the corresponding
property for the pair $(y,x)$, so we can consider only unordered pairs
$\{x,y\}$, and this reduces the number of intersection properties to
verify by a factor of 2. Still, the total number of intersection
properties is quadratic on the number of vertices 
and exponential
on the number of colors. 
This number gets too big
too quickly; however, many of these properties are redundant,
either because they are true for any group (for example, the
intersection of a group and one of its subgroups is the smaller
subgroup) or because they are a consequence of other intersection
properties.

Fortunately we may reduce the number of intersection properties to
check by following the same proof but using the {\em weak} path
intersection property instead of the strong one. Doing this we get
the following refinement of the previous theorem.

\begin{teo}\label{t:IntProp}
  Given a premaniplex $X$ and a voltage assignment $\xi$ such that
  $X^\xi$ is a maniplex, $X^\xi$ is the flag graph of a polytope if
  and only if
  \begin{equation}\label{eq:IntPropIntervalos}
    \xi(\Pi^{x,y}_{[0,m]}(X)) \cap \xi(\Pi^{x,y}_{[k,n-1]}(X)) =
    \xi(\Pi^{x,y}_{[k,m]}(X)),
  \end{equation}
  for all $k,m \in \{0,\ldots,n-1\}$ and all $x,y\in \F$.
\end{teo}

With Theorem~\ref{t:IntProp} the number of intersection properties to
check is now quadratic on the number of vertices and also quadratic on
the rank. We can still refine this theorem a little more with the following observations:
\begin{itemize}
    \item The cases $k=0$ and $m=n-1$ say that the intersection of the whole voltage group with some subset is the subset itself, so they are trivially true for every voltage assignment.

    \item 
    In Theorem~\ref{t:IntProp} one has to consider the case when $k>m$.
    In such case, $\xi(\Pi^{x,y}_{[k,m]}(X))$ is the trivial group when $x=y$ and the empty set when $x \neq y$.
    However, if the intersection property holds for $k=m+1$, that is, if $\xi(\Pi^{x,y}_{[0,m]}(X)) \cap \xi(\Pi^{x,y}_{[m+1,n-1]}(X)) =
    \xi(\Pi^{x,y}_{\emptyset}(X))$, then it also holds for $k>m+1$. So one may only verify the intersection property for $k\leq m+1$.

    \item If $y$ and $y'$ are in the same connected component of $X_{[k,m]}$ and ~(\ref{eq:IntPropIntervalos}) is satisfied for the pair $(x,y)$, then it is also satisfied for the pair $(x,y')$.
    To see this, let $W$ be a path from $y$ to $y'$ with colors in $[k,m]$ and notice that $\Pi^{x,y'}_I(X) = \Pi^{x,y}_I(X)W$ whenever $I$ contains $[k,m]$. By taking voltages we get that $\xi(\Pi^{x,y'}_I(X)) = \xi(W)\xi(\Pi^{x,y}_I(X))$.
    This means thatwe can get the intersection property for the pair $(x,y')$ by multiplying the one for pair $(x,y)$ by $\xi(W)$ on the left.
    So for each pair of numbers $(k,m)$ we only need to verify one intersection property for each pair of connected components of $X_{[k,m]}$.
\end{itemize}

Taking the previous observations into consideration, the maximum amount of necessary intersection properties to check is $v(v+1)/2\sum_{m=0}^{n-2} (m+1) =  v(v+1)n(n-1)/4$.
But to reach this bound we need for every vertex to be in a different connected component of $X_{[k,m]}$ for every pair $(k,m)$ with $0<k\leq m+1$, which is only possible when $v=1$. In other words, this bound is not tight for symmetry type graphs with more than one vertex.

\ifEjemploPropsTresOrb
Using Theorem~\ref{t:IntProp} on the premaniplex $3^1$ we get the following
intersection properties (in parenthesis the values of $k$ and $m$ and
a choice endpoints for each property):

\begin{eqnarray}\label{IntProp3^1s}
    \begin{aligned}
        \gen{\alpha_1, \beta_0} \cap \gen{\alpha_1,\beta_2^{\alpha_{02}}} &=& \gen{\alpha_1}\quad (k=1, m=1, x\to x)\nonumber\\
        \gen{\alpha_1, \beta_0} \cap \gen{\alpha_1^{\alpha_{02}}, \beta_2} &=& 1\quad (k=1, m=1, y\to y)\nonumber\\
        \gen{\alpha_1,\beta_0} \cap \alpha_{02}\gen{\alpha_1,\beta_2^{\alpha_{02}}} &=& \emptyset\quad  (k=1, m=1, x\to y)\nonumber\\
        \gen{\alpha_1,\beta_0} \cap \{\alpha_{02}\} &=& \emptyset\quad (k=2, m=1, x\to y)\nonumber\\
        \gen{\beta_0} \cap \gen{\alpha_1^{\alpha_{02}},\beta_2} &=& 1\quad (k=1, m=0, z\to z).
  \end{aligned}
\end{eqnarray}\elias{arreglar para que quede bien la referencia}
We have omitted the properties that are redundant or that hold trivially for any voltage assignment.

Applying Theorem~\ref{t:IntProp} to the premaniplex $3^1$ and using  Corollary~\ref{c:Mani3^1} we get the following result:

\begin{teo}\label{t:3^1}
  A group $\Gamma$ acts on a polytope with symmetry type graph $3^1$ if and only if it is generated by 4 involutions $\alpha_1,\alpha_{02}, \beta_0$ and $\beta_2$ such that $\beta_0$ commutes with $\beta_2$ and the intersection properties~(\ref{IntProp:3^1s} hold.
\end{teo}

A previous characterization for the automorphism groups of polytopes with symmetry type graph $3^1$ was given in~\missingref[el de 3 óbitas], however in that paper 7 intersection properties need to be checked, in contrast to the 5 presented here.
We can get the 7 conditions in~\missingref[el de 3 óbitas] if we use Theorem~\ref{t:IntProp-w} instead and make some choices for $I$
 and $J$. In this case, using Theorem~\ref{t:IntProp} only saves us 2 intersection properties, but in higher rank the difference would be way more significant.
\fi

\subsection{Constructing a polytope from the voltage
  group}\label{s:Construction}

We have seen how to recover a polytope from its flag graph (see
Theorem~\ref{t:PolyMani}) and when
$X^\xi$ is the flag graph of a polytope for a given premaniplex $X$
and a voltage assignment $\xi$ (see Theorem~\ref{t:IntProp}). By
concatenating the construction of $X^\xi$ from $X$ and $\xi$, and the
construction of a polytope $\p$ from $X^\xi$ we get a construction of
a polytope from $X$ and $\xi$.
In this section, we translate this to give a
construction only in terms of subgroups of $\Gamma$ and their
cosets. This will give the proof of the Main Theorem.

Let $f$ be a choice function on the connected subgraphs of $X$, that is, a function that assigns a base vertex to each such subgraph.
Let $C$
be a connected component of $X_I$ for some
$I\subset \{0,1,\ldots,n-1\}$ and let $x=f(C)$.
Let $\ol{C}:=(X^\xi)_I(x,1)$, that is, the connected
component of $(X^\xi)_I$ containing $(x,1)$.
Thus, if
$(x,\gamma) \in \ol{C}$, then there is a path $\os{W}$ from $(x,1)$ to
$(x,\gamma)$ which uses only colors in $I$, which implies that its projection is a
closed path $W$ in $\Pi^x_I(X)$ with voltage $\gamma$. This means that
when considering the action of $\Gamma$ on $X^\xi$, the stabilizer of
$\ol{C}$ coincides with $\xi(\Pi^x_I(X))$. If we now consider a coset
$\xi(\Pi^x_I(X)) \sigma$, this would be the set of elements of
$\Gamma$ that map $\ol{C}$ to $(X^\xi)_I(x,\sigma)$.

We know that the $i$-faces of the polytope that has $X^\xi$ as its
flag graph correspond to the connected components of
$(X^\xi)_{\ol{i}}$.
This makes natural the following construction:

Given $X$ and $\xi$ satisfying Theorem~\ref{t:IntProp} and 
a choice function $f$ on the connected subgraphs of $X$,
we construct a
partially ordered set $\p(X,\xi)$ with a rank function whose elements
of rank $i$ are the right cosets of groups of the type
$\xi(\Pi^x_{\ol{i}}(X))$ where
\[
x \in \{f(C)|C\ {\rm is\ a\ connected\
  component\ of}\ X_{\ol{i}}\}.
 \]
  We have to consider these groups as formal copies, that is, if $x$ and $x'$ are on different connected components we consider $\xi(\Pi^x_{\ol{i}}(X))$ to be different than $\xi(\Pi^{x'}_{\ol{i}}(X))$, even if, as groups, they might be equal.
  Likewise, if $j\neq i$ but $\xi(\Pi^x_{\ol{i}}(X))$ coincides with $\xi(\Pi^{x'}_{\ol{j}}(X))$ we consider them to be different in $\p(X,\xi)$.
  We could formalize this by saying that the elements of rank $i$ in $\p(X,\xi)$ are pairs $(C,\xi(\Pi^x_{\ol{i}}\gamma))$ with $C$ a connected component of $X_{\ol{i}}$, $x=f(C)$ and $\gamma \in \Gamma$, but since $x$ and $i$ already appear in the notation, we may assume that $\xi(\Pi^x_{\ol{i}}(X))\gamma$ stands for the pair $(X_{\ol{i}}(x),\xi(\Pi^x_{\ol{i}}(X))\gamma)$.

%

Now we shall define the order on $\p(X,\xi)$ as follows.
First, for all $i\in\{0,1,\ldots,n-1\}$ and every vertex $y$
in $X$ we look at the connected component $X_{\ol{i}}(y)$ and fix a
path going from its base vertex $x=f(X_{\ol{i}}(y))$ to $y$. We call
this path $W^y_i$ and we denote its voltage by
$\alpha_i^y:=\xi(W^y_i)$.

\begin{Def}[Order in $\p(X,\xi)$]\label{d:p(X,xi)}
  Given $0\leq i<j \leq n-1$, let $C$ be a
  connected component of $X_{\ol{i}}$ and $C'$ be a connected
  component of $X_{\ol{j}}$, and let $x$ and $x'$ be their respective
  base flags.  Given $\gamma, \gamma' \in \Gamma$, we say that
  \begin{eqnarray}
  \xi(\Pi^x_{\ol{i}}(X))\gamma < \xi(\Pi^{x'}_{\ol{j}}(X))\gamma'\ \mathrm{if} \ \mathrm{
  and}\ \mathrm{ only}\ \mathrm{ if} \
 \alpha^y_i \xi(\Pi^x_{\ol{i}}(X))\gamma \cap
  \alpha^y_j\xi(\Pi^{x'}_{\ol{j}}(X))\gamma' \neq \emptyset, \nonumber
    \end{eqnarray}
    for some $y \in C \cap C'$.
\end{Def}

\begin{teo}\label{t:Recons}
  Let $X$ be a premaniplex and $\xi:\Pi(X)\to \Gamma$ a voltage
  assignment satisfying Theorem~\ref{t:IntProp}. Let
  \[
  \p(X,\xi) := \{\xi(\Pi^x(C))\tau:C {\rm \ is\ a\ connected\ component
   \ of\ }
  X_{\ol{i}}, x=f(C), \tau \in \Gamma\},
  \]
  together with the order given in Definition~\ref{d:p(X,xi)}. Then
  $\p(X,\xi)$ is a polytope in which $\Gamma$ acts with symmetry type
  graph $X$.
\end{teo}

\begin{dem}


%

Let us denote $\p(X^\xi)$ as $\p$ and $\p(X,\xi)$ as $\Q$.
  Note that Theorem~\ref{t:IntProp} implies that $\Q$ is a polytope, and we want to show that $\p$ is also a polytope (such that $\p / \Gamma \cong X$).

  By the discussion in Section~\ref{s:VoltMani},  $\T(\p,\Gamma)\cong X$, since, by construction of $X^\xi$, we have that $X^\xi/\Gamma \cong X$.
 So in order to settle the theorem it is enough to find a poset isomorphism
  $\varphi:\Q\to \p$ such that it commutes with the action of
  $\Gamma$, that is, such that $\os{C}\sigma\varphi = \os{C}\varphi\sigma$ for all
  faces $\os{C}$ in $\Q$ and all $\sigma\in\Gamma$.

  Let $\os{C}$ be a face of $\p$. Hence,
  $\os{C}$ is a connected component of $X^\xi_{\ol{i}}$ for some color
  $i$. Let $C:=p(\os{C})$ and let $x:=f(C)$. This implies that $\os{C}$ has a flag
  of the type $(x,\gamma)$ for some $\gamma \in \Gamma$. Let $\os{K}$ be
  the connected component of $X^\xi_{\ol{i}}$
  that contains $(x,1)$.
Then, the set of elements of $\Gamma$
  that map $\os{K}$ to $\os{C}$  is
  $\xi(\Pi^x_{\ol{i}}(X))\gamma$. We want to identify $\os{C}$ with this
  coset, so we define $\os{C}\varphi:=\xi(\Pi^x_{\ol{i}}(X))\gamma$. Note that $\varphi:\Q\to\p$ is well defined, since if $(x,\gamma')$ is in $\os{C}$ then $\gamma'\gamma^{-1}$ stabilizes $\os{K}$, implying that $\gamma'\gamma^{-1} \in \xi(\Pi^x_{\ol{i}}(X))$.
  We want to prove that $\varphi$ is a poset isomorphism and that it
  commutes with $\Gamma$.

  Let us show first that $\varphi$ commutes with the action of
  $\Gamma$.
  Let $\os{C}$ be a face in $X^\xi$ and let $\sigma\in
  \Gamma$. By the definition of $\varphi$ we know that $\os{C}\varphi
  = \xi(\Pi^x_{\ol{i}}(X))\gamma$ where $i$ is the rank of $\os{C}$,
  $x=f(p(\os{C}))$ and $\gamma\in\Gamma$ is any element such that
  $(x,\gamma)\in \os{C}$. On the other hand $(\os{C}\sigma)\varphi =
  \xi(\Pi^{x'}_{\ol{i}}(X))\gamma'$ where $x'=f(p(\os{C}\sigma))$
  and $(x',\gamma')\in \os{C}\sigma$. Also note that since $\os{C}$
  and $\os{C}\sigma$ are in the same orbit, then
  $p(\os{C})=p(\os{C}\sigma)$, and thus $x=x'$. Furthermore, the flag
  $(x,\gamma\sigma) = (x,\gamma)\sigma$ is in $\os{C}\sigma$. This
  proves that $(\os{C}\sigma)\varphi =
  \xi(\Pi^x_{\ol{i}}(X))\gamma\sigma = (\os{C}\varphi)\sigma$.

  Now let us prove that $\varphi$ is an isomorphism of posets. Let
  $\os{C}$ and $\os{C}'$ be incident faces of $\p$ of ranks $i$ and
  $j$, respectively, with $i<j$ (therefore,
  $\os{C}<\os{C}'$). Hence, there is a flag $(y,\tau)$ in $\os{C} \cap
  \os{C}'$, which in turn implies that its first entry, $y$, must be in $C \cap C'$ where
  $C=p(\os{C})$ and $C'=p(\os{C}')$.

  Note that the path $W^y_i$ is contained in $C$ while
  $W^y_j$ is contained in $C'$. Then, these paths have lifts
  $\os{W}^y_i$ and $\os{W}^y_j$ respectively, that go from
  $(x,(\alpha^y_i)^{-1}\tau)$ and $(x',(\alpha^y_j)^{-1}\tau)$,
  respectively, to $(y,\tau)$. Observe that $\os{W}^y_i$ is contained in
  $\os{C}$ and $\os{W}^y_j$
  is contained in $\os{C}'$. Thus, $(\alpha^y_i)^{-1}\tau \in
  \Pi^x_{\ol{i}}(X)\gamma$ and $(\alpha^y_j)^{-1}\tau \in
  \Pi^{x'}_{\ol{j}}(X)\gamma'$. Therefore
  \[
    \tau \in \alpha^y_i \xi(\Pi^x_{\ol{i}}(X))\gamma \cap
    \alpha^y_j \xi(\Pi^{x'}_{\ol{j}}(X))\gamma' = \alpha^y_i
    (\os{C}\varphi) \cap \alpha^y_j (\os{C}'\varphi);
  \]
  but this means that $\os{C}\varphi<\os{C}'\varphi$ in $\Q$.

  Conversely, suppose that $\os{C}\varphi<\os{C}'\varphi$ in
  $\Q$. We want to
  show that $\os{C}<\os{C}'$ in $X^\xi$. Let us write $\os{C}\varphi =
  \xi(\Pi^x_{\ol{i}}(X))\gamma$, where $i$ is the rank of $\os{C}$, $x:=f(p(\os{C}))$ and $\gamma$ is an element of
  the voltage group such that $(x,\gamma)\in \os{C}$. Similarly, we
  write $\os{C}'\varphi =
  \xi(\Pi^{x'}_{\ol{j}}(X))\gamma'$, where $j$ is the rank of $\os{C}'$,
  $x:=f(p(\os{C}'))$ and $\gamma'$ is an
  element of the voltage group such that $(x',\gamma')\in \os{C}'$.

  By hypothesis $\alpha^y_i \xi(\Pi^x_{\ol{i}}(X))\gamma$ and $\alpha^y_j
  \xi(\Pi^{x'}_{\ol{j}}(X))\gamma'$ have non-empty intersection for some
  $y\in C\cap C'$. Let $\tau$
  be an element in such intersection. Then $(\alpha^y_i)^{-1}\tau \in
  \xi(\Pi^x_{\ol{i}}(X))\gamma$. This implies that
  $(x,(\alpha^y_i)^{-1}\tau)$ is in the same connected component of
  $(X^\xi)_{\ol{i}}$ as $(x,\gamma)$, that is $(x,(\alpha^y_i)^{-1}\tau)
  \in \os{C}$. But at the same time, there is a
  lift of $W^y_i$ that connects $(x,(\alpha^y_j)^{-1}\tau)$ with
  $(y,\tau)$, and since $W^y_i$ does not use the color $i$, its lift is
  contained in $\os{C}$, which proves that
  $(y,\tau)\in\os{C}$. Similarly, the fact that
  $(\alpha^y_j)^{-1}\tau \in \xi(\Pi^{x'}_{\ol{j}}(X))\gamma$ implies
  that $(y,\tau)\in \os{C}'$. Thus, we have proved that $\os{C}\cap
  \os{C}'$ is not empty, or in other words $\os{C}<\os{C}'$ in
  $\p$.

  Therefore, $\varphi$ is an isomorphism and the theorem follows. \qed

\end{dem}

\section{Example: Caterpillars}\label{s:Caterpillars}

If one wants to build polytopes from premaniplexes in the way we have described in this paper, it is natural to start with infinite families of premaniplexes.
One could, of course, then start with premaniplexes with a fixed number $k$ of vertices. For example, one may want to construct 2-orbit polytopes (see \cite{IsaEgon2orbit}).
If, on the other hand, we do not want to limit the number of vertices of the premaniplexes in the family a first step could be to consider trees. However, since semi-edges are considered as cycles, there are no premaniplexes whose underlying graph (without colours) is a tree. Thus, we study the closest thing to them: premaniplexes that are trees with an unlimited number of semi-edges.
Hence, we define a {\em caterpillar} as a premaniplex in which every cycle  
is a semi-edge.
In other words, a caterpillar is a premaniplex $X$ with a
unique spanning tree.

In particular, a
caterpillar does not have pairs of parallel links (edges joining different vertices). If there are three
links  incident to one vertex, at
least two of them must have colors differing by more than 1, which
would imply that there is a 4-cycle. This implies that caterpillars
consist in fact of a single path $P$ (which we will call {\em the
  underlying path of $X$}) and lots of semi-edges. Of course, the
colors of two consecutive edges on the path must differ by exactly
one, otherwise there would be a 4-cycle.
We note here that the term {\em caterpillar} has been used in the graph theory literature for a very similar but slightly different concept.

Throughout this section, unless otherwise stated, $X$ will denote a finite caterpillar (that is, one with a finite number of vertices) with underlying path $P$, and its vertices will be labeled by $x_0, x_1,\ldots, x_k$, ordered as they are visited by $P$.
Furthermore, we denote by $(x_i,j)$ the dart of color $j$ at $x_i$ and by $c_i$ the color of the link connecting $x_{i-1}$ and $x_i$.

\begin{figure}
  \begin{center}
    \includegraphics[width=5cm]{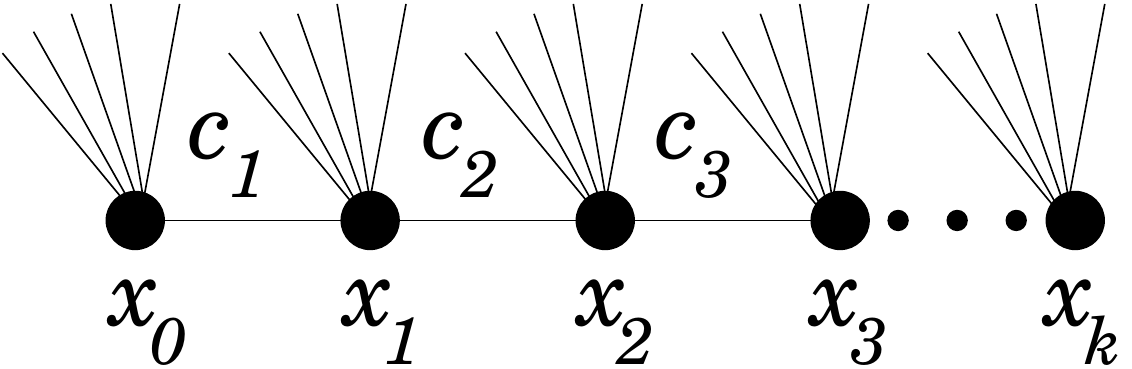}
    \caption{\label{f:caterpillar} A finite caterpillar.}
  \end{center}
\end{figure}

\subsection{Caterpillar coverings}

We want to construct polytopes from caterpillars.
As often when constructing maniplexes and polytopes via voltage assignments from a premaniplex $X$, the derived maniplex might not have $X$ as the symmetry type graph with respect to the full automorphism group.
However, the actual symmetry type graph is a quotient of $X$.
For this reason, in this section we study quotients of caterpillars.

Given a caterpillar, let us call an {\em endpoint} a vertex incident to just one link (that is, an endpoint of the underlying path).
Note that a caterpillar is finite if and only if it has exactly
two endpoints.
Every symmetry of the caterpillar must map endpoints to endpoints. If a caterpillar is finite there is at most one non-trivial symmetry and its action on the vertices $x_0,x_1,\ldots,x_k$ is given by $x_j\mapsto x_{k-j}$.
We call a finite caterpillar {\em symmetric} if it has a non-trivial symmetry.

A {\em word} $w$ in $\{0,1,\ldots,n-1\}$ is simply a finite sequence $w=a_1a_2\ldots a_t$, with $a_i\in \{0,1,\ldots,n-1\}$ for each $i=1,2,\ldots,t$.
The {\em inverse} of a word $w$ is the word $w^{-1}$ that has the same colors as $w$ but written in reverse order.
That is, if $w=a_1a_2\ldots a_t$ then $w^{-1} = a_ta_{t-1}\ldots a_1$. A word is said to be {\em reduced} if it has no occurrence of the same color twice in a row; in other words, $w=a_1a_2\ldots a_t$ is reduced if and only if $a_{i+1}\neq a_i$ for all $i=1,2,\ldots,t-1$.
We shall work with reduced words from now on.
A word $w=a_1a_2\ldots a_t$ is a {\em palindrome} if  $a_i = a_{t+1-i}$ for all $i\in \{1,2,\ldots,t\}$.
A palindrome word of even length can
be written as $vv^{-1}$ for some word $v$ and cannot be a
reduced word.
A palindrome word of odd length can always be written as
$w=vav^{-1}$ for some color $a$ and some word $v$.

Given a segment $[x,y]$ in a caterpillar, its {\em underlying word} is
the word $w$ consisting of the colors of the links in the path that
goes from $x$ to $y$. When we speak of the {\em underlying word} of a
caterpillar $X$ we are referring to the underlying word of its
underlying path in a fixed orientation.
Hence, we say that a segment $[x,y]$ is a {\em palindrome} if its underlying word
$v$ is a palindrome.

\begin{prop}\label{p:CubsCat}
  Let $X$ be a finite caterpillar and let $Y$ be a premaniplex not
  isomorphic to $X$ 
  such that there is a premaniplex homomorphism 
  $h:X\to Y$. Then $Y$ is a caterpillar. Moreover, if $Y$ has at least
  2 vertices and $S=c_1c_2\ldots c_k$ is the underlying word of $X$,
  then there is some $r<k$ such that $w=c_1c_2\ldots c_r$ is the
  underlying word of $Y$ and one of the following statements is true:
  \begin{enumerate}
  \item There exist colors $a_1,a_2,\ldots,a_t \in \{c_r+1,c_r-1\}$ and
    $b_1,b_2\ldots,b_{t-1} \in \{c_1+1,c_1-1\}$ such that
    $S=wa_1w^{-1}b_1wa_2w^{-1}b_2\ldots b_{t-1}wa_tw^{-1}$.
  \item There exist colors $a_1,a_2,\ldots,a_t\in \{c_r+1,c_r-1\}$ and
    $b_1,b_2\ldots,b_t \in \{c_1+1,c_1-1\}$ such that
    $S=wa_1w^{-1}b_1wa_2w^{-1}b_2\ldots b_{t-1}wa_tw^{-1}b_tw$.
  \end{enumerate}

  In any case, if $i \equiv j \md 2r+2$ then $h(x_i) = h(x_j)$. Also
  if $i \equiv -j-1 \md 2r+2$ 
  then $h(x_i) = h(x_j)$
  .
\end{prop}

\begin{figure}
  \centering
  \begin{subfigure}[t]{0.45\linewidth}
    \includegraphics[width=6cm]{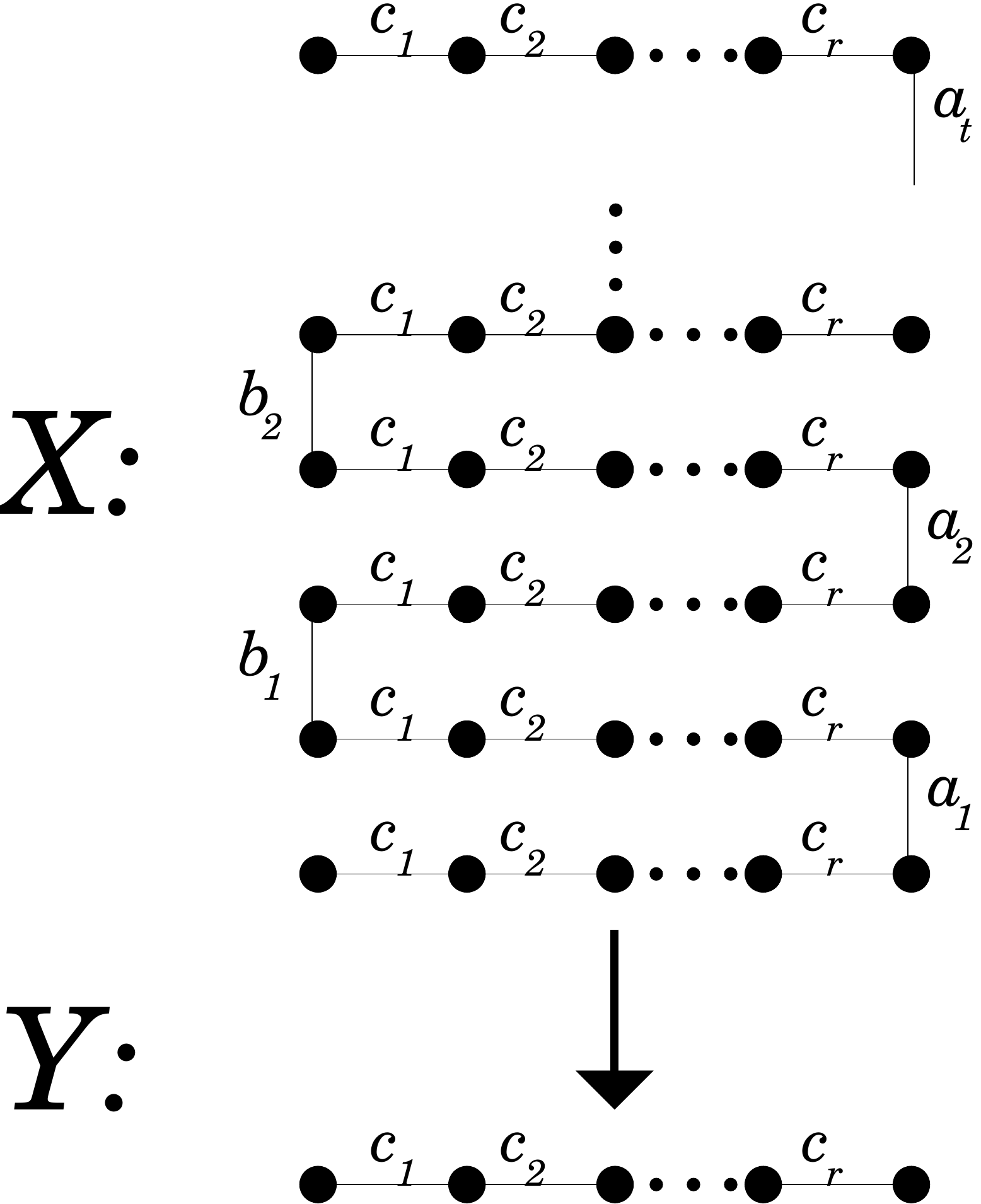}
    \caption{Case 1}
  \end{subfigure}
  \hfill
  \begin{subfigure}[t]{0.45\linewidth}
    \includegraphics[width=6cm]{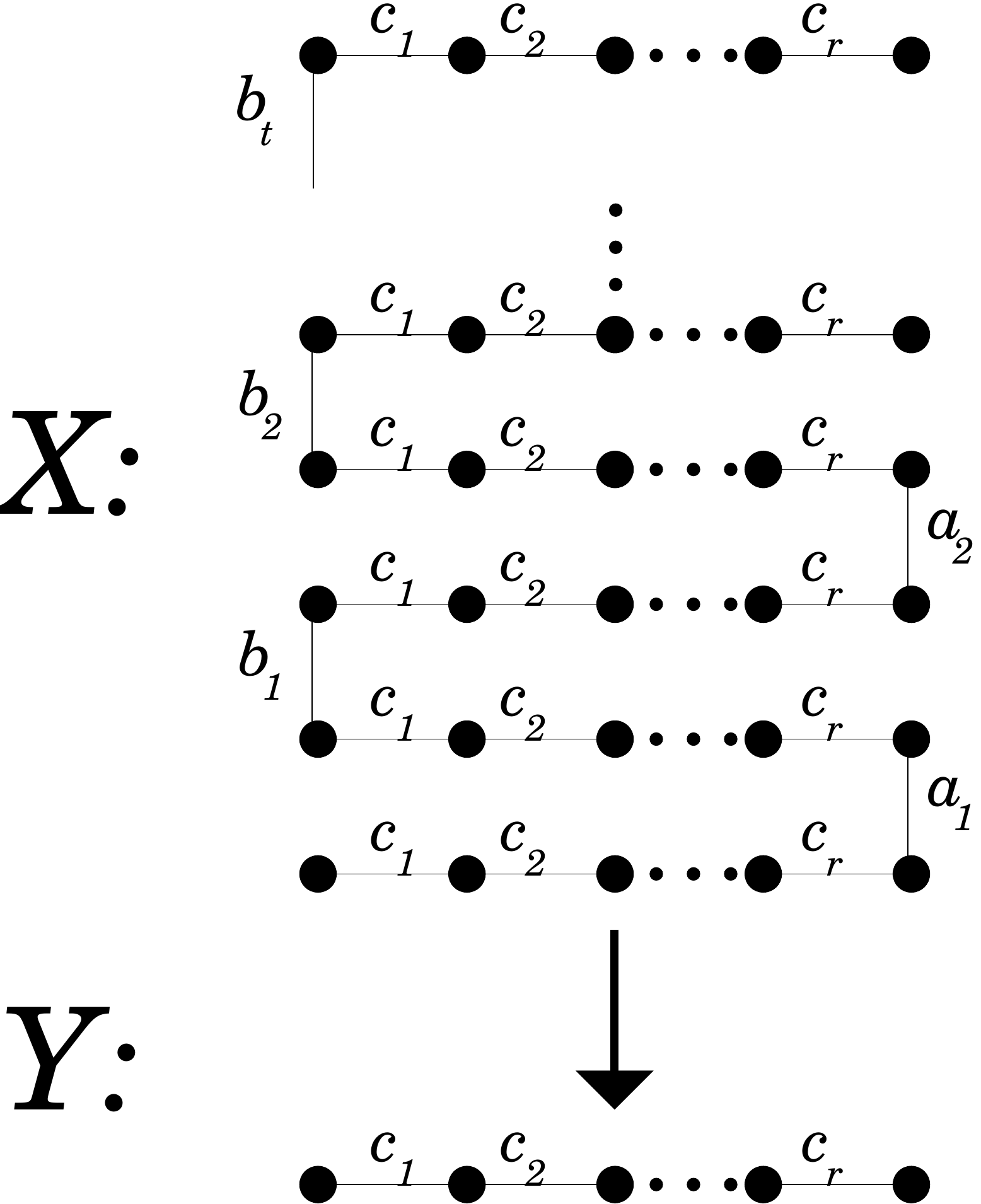}
    \caption{Case 2}
  \end{subfigure}
  \caption{The caterpillar $X$ covers the caterpillar $Y$ if and only
    if it ``folds'' into it.}
  \label{f:FoldedCat}
\end{figure}

Before proving Proposition~\ref{p:CubsCat} let us remark that it
simply means that the quotients of a caterpillar $X$ are those
caterpillars $Y$ such that $X$ can be ``folded'' into $Y$. We
illustrate this concept in Figure~\ref{f:FoldedCat}: the semi-edges
are not drawn and the names of the vertices have been omitted, but the
idea is that $X$ must be ``folded'' into ``layers'' of $r+1$ vertices
and then each vertex will be projected to the vertex on $Y$ in the
same horizontal coordinate. The layer $\ell$ consists of the vertices
$x_i$ where $\floor{\frac{i}{r+1}}=\ell$ ($\floor{x}$ denotes the
integer part of $x$). Even layers go from left to right, while odd
layers go from right to left, hence the underlying word of even layers
is $w=c_1c_2\ldots c_r$ while the underlying word of odd layers is
$w^{-1}=c_rc_{r-1}\ldots c_1$.

Now let us proceed with the proof.

\begin{dem} {\bf of Proposition~\ref{p:CubsCat}}
  Let us consider the equivalence relation $\sim_h$ on $X$ given by  $x\sim_hy$ if and only if $h(x)=h(y)$.
  Recall that by the definition of a premaniplex homomorphism $x\sim_h y$ implies that $x^m\sim_h y^m$ for every monodromy $m$.
  Therefore $Y$ is isomorphic to $X/\sim_h$.

  We already know that all premaniplex homomorphisms are surjective
  (see Section~\ref{sec:Polytopes-and-maniplexes}). By hypothesis $h$ is not an isomorphism,
  which implies it cannot be injective. Let $x$ and $y$ be two
  different vertices on $X$ such that $x\sim_h y$. Let
  $x_0,x_1,\ldots, x_k$ be the sequence of vertices in the underlying
  path of $X$. There is a monodromy $m$ such that $x^m=x_0$, and so
  $x_0\sim_h y^m$ and $y^m$ is different than $x_0$. Now let $q$ be
  the minimum positive number such that $x_q\sim_h x_0$. We know that
  $q$ exists because $x_0\sim_h y^m$. Note that if
  $q=1$ we would have that $x_0^i\sim_h x_0$ for all $i$, which would
  imply that the same is true for $x_1$ and in turn also for
  $x_2=x_1^{c_2}$ and so on. This would mean that all the vertices of
  $X$ are equivalent, meaning that $Y$ has only one vertex and the
  proposition follows. So we may assume that $q>1$.

  We know that $x_{q-1} = x_q^{c_q}$ is equivalent $x_0^{c_q}$. If
  $c_q \neq c_1$, we would have that $x_0^{c_q}=x_0$. This would imply
  that $x_{q-1}$ is equivalent to $x_0$, contradicting the minimality
  of $q$. So we have proven that $c_q = c_1$.

  Now note that if for some $\ell$ we have that $x_\ell\sim x_1$, then
  $x_\ell^{c_1}\sim x_1^{c_1}=x_0$. In particular this tells us that
  for $\ell<q-1$, $x_\ell$ cannot be equivalent to $x_1$. Now we can
  use the same argument we used to prove that $c_q=c_1$ to prove that
  $c_{q-1}=c_2$. Analogously we can prove that
  $c_3 = c_{q-2}, c_4=c_{q-3}$ and so on. In other words, $[x_0,x_q]$
  is a palindrome. Since for all $i$, $c_i$ and $c_{i+1}$ are
  different, then $q$ is odd, say $q=2r+1$. Let $v$ be the underlying
  word of the segment $[x_0,x_q]$. Then $v$ may be written as
  $v=wc_{r+1}w^{-1}$ for some word $w = c_1c_2\ldots c_r$. Let us call
  $a_1:=c_{r+1}$ and note that $a_1\in \{c_r+1,c_r-1\}$.

  For all $i=0,1,\ldots, k$, let us denote by $\widehat{i}$ the
  residue of dividing $i$ by $2r+2$. We will prove by induction on $i$
  that $x_i\sim_h x_{\widehat{i}}$ for all $i=0,1,\ldots, k$ and that
  $c_{i} = c_{\widehat{i}}$ when $i$ is not divisible by $r+1$,
  $c_i \in \{c_1+1,c_1-1\}$ if $i$ is an odd multiple of $r+1$ and
  $c_i= \{c_r+1,c_r-1\}$ when $i$ is an even multiple of $r+1$.

  Let our induction hypothesis be that $x_{\ell}\sim_h
  x_{\widehat{\ell}}$ for all $\ell<i$, and that
  $c_{\ell} = c_{\widehat{\ell}}$ if $\ell$ is not divisible by $r+1$.


  Let us start with the case when $i\equiv 0 \md 2r+2$. In this case
  we want to prove that $x_i \sim_h x_{\widehat{i}}$ and that
  $c_i\in \{c_1+1,c_1-1\}$. By our induction hypothesis we know that
  $x_{i-1}\sim_hx_{2r+1} \sim_h x_0$ and $c_{i-1}=c_{2r+1}=c_1$. In
  particular $c_i\in \{c_1+1,c_1-1\}$.  This implies that
  $x_i=x_{i-1}^{c_i}\sim_h x_0^{c_i}=x_0=x_{\,\widehat{i}}$ (since
  $c_i\neq c_1$). Thus $i$ satisfies our claim.

  Now let us proceed with the case when $i$ is an odd multiple of
  $r+1$, that is $\widehat{i} = r+1$. In this case we want to prove
  that $x_i \sim_h x_{\widehat{i}}$ and that $c_i\in
  \{c_r+1,c_r-1\}$. Our induction hypothesis tells us that
  $x_{i-1}\sim_h x_r$ and that $c_{i-1}=c_r$. Hence
  $c_i \in \{c_r+1,c_r-1\}$. Note that one of the colors in
  $\{c_r+1,c_r-1\}$ is actually $c_{r+1}$, while the other is the
  color of a semi-edge incident to $x_r$. Since $x_{r+1}\sim_h x_r$ we
  have that
  $x_i=x_{i-1}^{c_i}\sim_h x_r^{c_i}\sim_h x_r\sim_h x_{r+1} =
  x_{\,\widehat{i}}$. Thus $i$ satisfies our claim.

  Finally let us prove our claim for the case when $i$ is not
  divisible by $r+1$. Our induction hypothesis tells us that
  $x_{i-1}\sim_h x_{\widehat{i-1}}$. Note that
  $\widehat{i} = \widehat{i-1}+1$. Since $i$ is not a multiple of
  $r+1$ we know that
  $x_{\,\widehat{i}}=x_{\widehat{i-1}}^{c_{\,\widehat{i}}}$ is not
  equivalent to $x_{\widehat{i-1}}$. This implies that
  $x_{i-1}^{c_{\,\widehat{i}}}$ is not equivalent to $x_{i-1}$, and
  since it is adjacent to $x_{i-1}$ it must be equal to either $x_i$
  or $x_{i-2}$. If $i\equiv 1 \md r+1$ then by induction hypothesis
  $x_{i-1}\sim_h x_0$, but recall that $x_0\sim_h x_{2r+1}=x_q$ by
  definition, and by induction hypothesis
  $x_{i-2}\sim_h x_{\widehat{i-2}} = x_{2r+1}$. This means that
  $x_{i-2}\sim_h x_{i-1}$, so $x_{i-1}^{c_{\,\widehat{i}}}$ must be
  $x_i$, implying that $ c_i=c_{\,\widehat{i}}$. If
  $i\not\equiv 1 \md r+1$ then our induction hypothesis tells us that
  $c_{i-1}=c_{\widehat{i-1}} \neq c_{\,\widehat{i}}$, and since
  $x_{i-2}=x_{i-1}^{c_{i-1}}=x_{i-1}^{c_{\widehat{i-1}}}$, so the only
  possibility is that $x_{i-1}^{c_{\,\widehat{i}}} = x_i$, implying that
  $c_i=c_{\,\widehat{i}}$ and that $x_i\sim_h x_{\,\widehat{i}}$. Thus,
  $i$ satisfies our claim. Note that we have also proven that if $i$
  is not divisible by $r+1$ then $x_{i-1}$ cannot be an endpoint of
  $X$.

  We have proved that $x_i\sim_h x_{\,\widehat{i}}$ for all
  $i=0,1,\ldots, k$. This implies automatically that if $i\equiv j \md
  2r+2$ then $x_i\sim_h x_j$. Moreover, if $i\equiv -j-1 \md 2r+2$,
  then $\widehat{i} = 2r+1-\widehat{j}$. Now since $v$ is a
  palindrome, we know that $x_\ell \sim_h x_{2r+1-\ell}$ for all
  $\ell=0,1,\ldots, 2r+1$, in particular $x_{\,\widehat{i}} =
  x_{2r+1-\widehat{j}} \sim_h x_{\,\widehat{j}}$. This, together with the
  fact that $x_i \sim_h x_{\,\widehat{i}}$ and $x_j \sim_h
  x_{\,\widehat{j}}$, implies that $x_i \sim_h x_j$.

  We have already proved that $S=wa_1w^{-1}b_1wa_2w^{-1}\ldots$ and it
  ends after an occurrence of $w$ or $w^{-1}$. To end the prove note
  that since $h$ is surjective we already know exactly what
  premaniplex $Y$ is: It is a caterpillar with vertex sequence
  $h(x_0), h(x_1),\ldots, h(x_r)$. In fact, for $i=1,\ldots,r-1$ and
  $j=0,1,\ldots n-1$ we know that $h(x_i)$ is different from
  $h(x_i)^j=h(x_i^j)$ if and only if $j\in\{c_{i-1},c_i\}$; and if
  $i=0$ (resp. $r$) then $h(x_i)$ is different from $h(x_i)^j$ if and
  only if $j=c_1$ (resp. $c_r$). \qed


\end{dem}





If we look closely at the proof of Proposition~\ref{p:CubsCat} we
notice that we have not actually used the fact that $X$ is finite, but
only the fact that it has at least one endpoint. So the proposition
may be generalized as follows:

\begin{prop}\label{p:CubsCat1end}
  Let $X$ be a caterpillar with at least one endpoint and let $Y$ be a premaniplex
  such that there is a premaniplex homomorphism (or covering) $h:X\to
  Y$. Then $Y$ is a caterpillar. Moreover, if $S=c_1c_2\ldots $ is
  the color sequence of $X$ starting at its endpoint, there is some
  $r$ such that $w=c_1c_2\ldots c_r$ is the underlying word of $Y$
  and there exist colors $a_1,a_2,\ldots \in \{c_r+1,c_r-1\}$ and
    $b_1,b_2\ldots \in \{c_1+1,c_1-1\}$ such that
    $S=wa_1w^{-1}b_1wa_2w^{-1}b_2\ldots $. If $X$ is finite it ends
    after an occurrence of either $w$ or $w^{-1}$. If $i \equiv j \md
    2r+2$ then $h(x_i) = h(x_j)$. Also if $i \equiv -j \md 2r+1$ and
    $\floor{\frac{i}{r+1}} \not\equiv \floor{\frac{j}{r+1}} \md 2$
    then $h(x_i) = h(x_j)$.
\end{prop}


Let us look now at the degenerate cases. If $Y$ is only one vertex we
can think that the previous theorems still hold but with $w$ being the
empty word, and since $c_1$ and $c_r$ do not exist, we get rid of the
restrictions $a_i\in \{c_r+1,c_r-1\}$ and $b_i\in \{c_1+1,c_1-1\}$. In
Figure~\ref{f:FoldedCat} $X$ would have only one vertex and no links
per layer; essentially $X$ would be ``standing up'' instead of being
folded. If $Y$ is isomorphic to $X$ then $S=w$ and we do not have any $a_i$
or $b_i$. In Figure~\ref{f:FoldedCat} $X$ would have only one layer.

If one would try to prove an analogous result as \cref{p:CubsCat1end} for caterpillars with no ends, the only additional case to consider would be that when $Y$ is not a caterpillar but a premaniplex consisting of a cycle with underlying word $w$ and all the edges not in the cycle are semiedges. In this case $X$ would have the underlying infinite word $\cdots w w w \cdots$.
\subsection{Caterpillars as STG of polytopes}

Given a caterpillar $X$, we want to assign voltages to the semi-edges of $X$ in order to get the flag graph of a polytope as the derived maniplex.
To this end, we note that in caterpillars the conditions of Theorems~\ref{t:IntProp-w} and~\ref{t:IntProp} can be simplified, as stated in the following lemma:

\begin{lema}\label{l:IntPropCat}
  Let $X$ be a caterpillar and let $\xi:\Pi(X)\to\Gamma$ be a voltage
  assignment such that all the darts in the underlying path of $X$
  have trivial voltage. Then the following statements are equivalent.
  \begin{enumerate}
  \item $X^\xi$ is polytopal.
  \item For every vertex $x$ in $x$ and all sets
    $I,J\subset\{0,1,\ldots,n-1\}$ the equation
    \[
      \xi(\Pi_I^x(X))\cap \xi(\Pi_J^x(X)) = \xi(\Pi^x_{I\cap J}(X))
    \]
    holds.
    \item For every vertex $x$ in $X$ and all
      $k,m\in\{0,1,\ldots,n-1\}$ the equation
      \[
        \xi(\Pi_{[0,m]}^x(X))\cap \xi(\Pi_{[k,n-1]}^x(X)) =
        \xi(\Pi^x_{[k,m]}(X))
      \]
      holds.
  \end{enumerate}
\end{lema}

\begin{dem}
  To prove this we have to see that for every set of colors
  $I\subset \{0,1,\ldots, n-1\}$ the set $\xi(\Pi^{x,y}_I(X))$ is either
  $\xi(\Pi^x_I(X))$ or empty, depending on whether or not the segment of
  $P$ that goes from $x$ to $y$ (which we will denote $[x,y]$) uses or
  not only colors in $I$. In fact if $[x,y]$ uses only colors in
  $I$ then
  \[
    \xi(\Pi^{x,y}_I(X)) = \xi(\Pi^x_I(X) [x,y]) =
    \xi([x,y])\xi(\Pi^x_I(X)) = 1\cdot \xi(\Pi^x_I(X)) =
    \xi(\Pi^x_I(X)).
  \]

  Note also that $\xi(\Pi^{x,y}_I(X)) \cap \xi(\Pi^{x,y}_J(X))$
  is empty if and only if one of the factors is empty. These
  observations together with Theorems~\ref{t:IntProp-w} prove the
  equivalence between conditions 1 and 2, and with
  Theorem~\ref{t:IntProp} we prove the equivalence between conditions
  1 and 3, thus proving the lemma. \qed

\end{dem}

Recall that a {\em Boolean group} (or an {\em elementary abelian $2$-group}) is a group in which every
non-trivial element has order exactly 2. Thus, all Boolean
groups are abelian, and finitely generated Boolean groups are
isomorphic to a direct product of cyclic groups of order 2.

\begin{prop}\label{p:Caterpillars1}
  Every caterpillar is the quotient of the flag graph of a polytope
  by a Boolean group.
\end{prop}

\begin{figure}
  \centering
  \includegraphics[width=7.5cm]{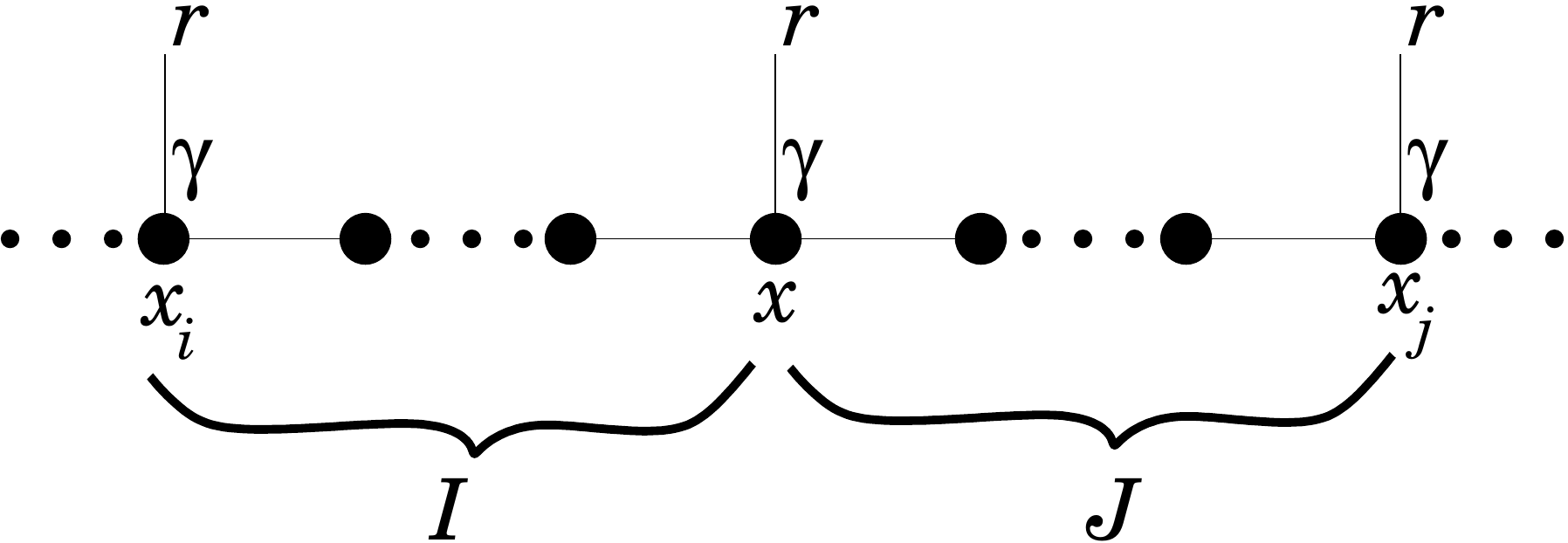}
  \caption{If $x\in[x_i,x_j]$ then $(x,r)$ has voltage
    $\gamma\in \xi(\Pi^x_{I\cap J}(X))$.}
  \label{f:IntPropCat_caso1}
\end{figure}

\begin{figure}
  \centering
  \includegraphics[width=7.5cm]{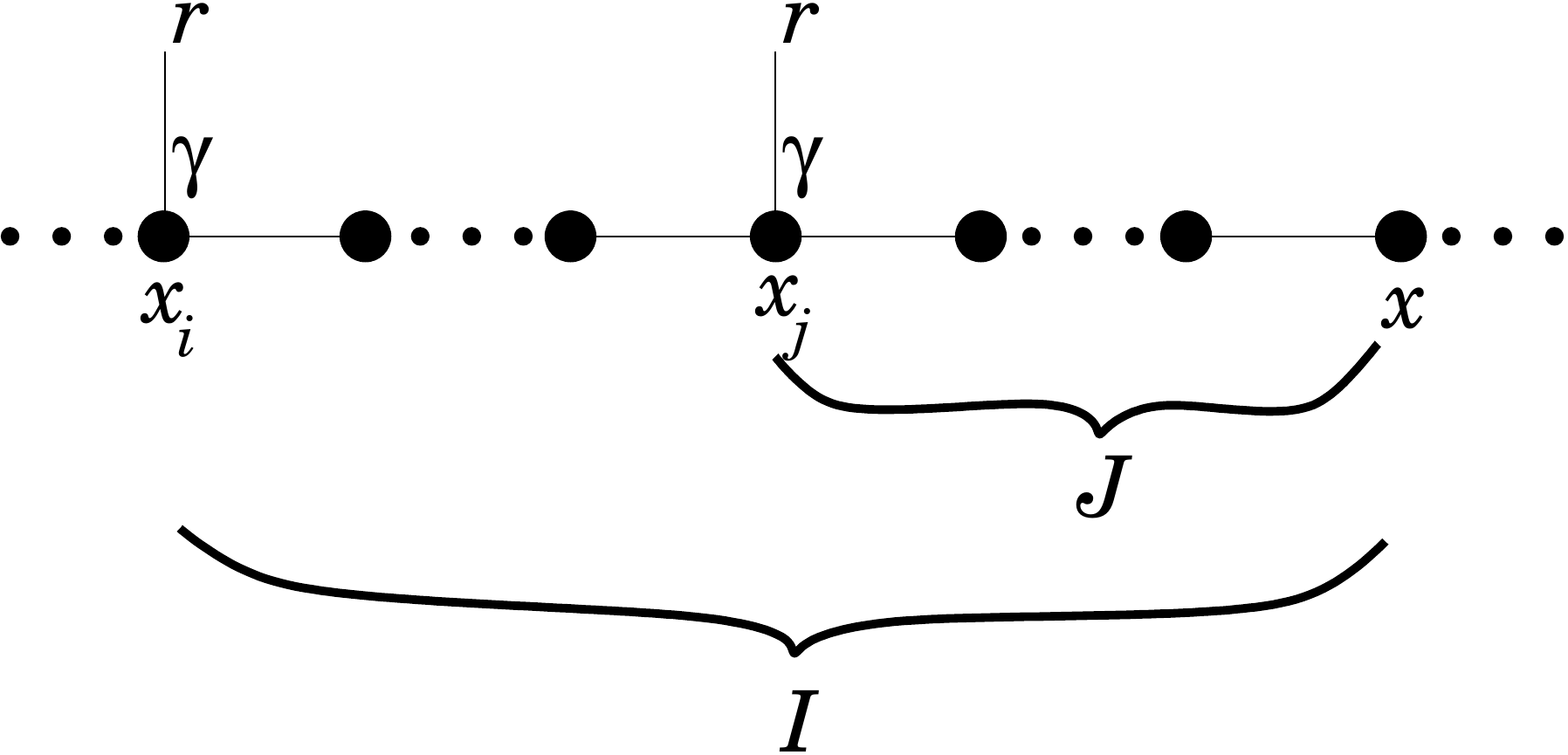}
  \caption{If $x_j\in[x_i,x]$ then $e'=(x_j,r)$ is in
    $X_{I\cap J}(x)$ and has voltage $\gamma$.}
  \label{f:IntPropCat_caso2}
\end{figure}

\begin{dem}
  Let $X$ be a caterpillar.
  We shall find a Boolean group $B$ and a voltage assignment $\xi:\Pi(X)\to B$ such that $X^\xi$ is polytopal in the following way.

  First assign trivial voltage to all links of $P$ and a different independent voltage to each semi-edge incident to $x_0$.
  We shall give voltage assignments to the semi-edges $(x_i,j)$, recursively on $i$.
  If $|j-c_i|>1$ we assign the same voltage to $(x_i,j)$ as we did to $(x_{i-1},j)$.
  On the ohter hand, if $|j-c_i|=1$, we assign a new element as the voltage $(x_i,j)$ independent from all the voltages of previous darts.

  \if
  To do so, define $\delta_i:=c_{i+1}-c_i$ for $i=1,\ldots,k-1$, and note that $\delta_i \in  \{-1,1\}$
  .
  There are two cases to consider to assign the voltage of $(x_i,j)$, for $i<k$:
  \begin{itemize}
      \item $j\neq c_i-\delta_i$.
      In this case $(x_{i-1},j)$ is a semi-edge because $j$ differs from $c_i$ in more than $1$, and the links incident to $x_{i-1}$ have colors $c_i$ and  $c_{i-1}=c_i+\delta_{i-1}$.
      This implies that the path $(x_{i-1},c_i)(x_i,j)(x_i,c_i)(x_{i-1},j)$ must have trivial voltage.
      Hence, we assign to the semi-edge $(x_i,j)$ the same voltage as the one that $(x_{i-1},j)$ has.
      \item $j= c_i-\delta_i$. In this case, we assign a new generator of a Boolean group as the voltage of  $(x_i,j)$.
  \end{itemize}
  Finally, since $\delta_k$ is undefined, we assign new
  generators as the voltages of the darts $(x_k,c_k+1)$ and
  $(x_k,c_k-1)$
  (if they exist) and  for the other semi-edges on $x_k$ copy the
  voltage of the semi-edge of the same color on $x_{k-1}$.
  \fi

  Let us call
  the resulting voltage group $B$ and denote this
  voltage assignment by $\xi$. By Lemma~\ref{l:VoltMani}, $X^\xi$ is a maniplex.

  Note further that $\xi$ satisfies that if $i\in
  \{1,2,\ldots,k\}$; $r,s\in \{0,1,\ldots,n-1\}$ and $(x_{i-1},r)$ and
  $(x_i,s)$ are semi-edges, then $\xi(x_i,s) = \xi(x_{i-1},r)$ if and
  only if $r=s$ and $|r-c_i|\neq 1$. It also satisfies that if
  $i<\ell<j$ and $(x_i,r)$ and $(x_j,r)$ are semi-edges such that
  $\xi(x_i,r) = \xi(x_j,r) = \gamma$, then $(x_\ell,r)$ is also a
  semi-edge and $\xi(x_\ell,r) = \gamma$.

  We claim that the second statement of Lemma~\ref{l:IntPropCat} is satisfied, and
  hence, $X^\xi$ is the flag graph of a polytope.
To show that this is the case, start by noticing that $\xi(\Pi^x_I)$ is the group generated by the voltages of the semi-edges on the component $X_I(x)$.

  Given $I, J \subseteq \{0,1, \dots, n-1\}$, suppose that for some vertex $x$ there is a semi-edge $e$ in
  $X_I(x)$ and a semi-edge $e'$ in $X_J(x)$ with
  $\xi(e) = \xi(e') = \gamma$ for some $\gamma \in \Gamma$. If $e=e'$
  then $e \in X_{I\cap J}(x)$. If $e\neq e'$ then we have that for
  some $i,j\in \{0,\ldots,k-1\}$ and some $r \in \{0,\ldots,n-1\}$
  occurs that $e=(x_i,r)$ and $e'=(x_j,r)$; in particular
  $r\in I\cap J$.
  If $x\in [x_i,x_j]$ then, as previously observed,
  $\xi(x,r) = \gamma$ (see Figure~\ref{f:IntPropCat_caso1}) and since
  $(x,r) \in X_{I \cap J}(x)$ this means $\gamma$ is a generator of
  $\xi(\Pi^x_{I\cap J})$. If $x\notin [x_i, x_j]$ consider without
  loss of generality that $x_i$ is further away from $x$ than $x_j$,
  in other words, that $x_j \in [x_i,x]$ (see Figure~\ref{f:IntPropCat_caso2}). Then
  $[x_i,x]$ uses only colors in $I$ and $[x_j,x] \subset [x_i,x]$ uses
  only colors in $J$. Then since $[x_j,x] \subset [x_i,x]$ we know
  that $[x_j,x]$ uses only colors in $I \cap J$, meaning that
  $e' \in X_{I \cap J}(x)$ and since its voltage is $\gamma$, this
  shows that $\gamma$ is always a generator of
  $\xi(\Pi^x_{I \cap J})$. Therefore, we have shown that if $\gamma$ is a
  generator of both $\xi(\Pi^x_I(X))$ and $\xi(\Pi^x_J(X))$, then it
  is also a generator of $\xi(\Pi^x_{I \cap J}(X))$.

  Now let $\sigma\in \xi(\Pi^x_I(X)) \cap \xi(\Pi^x_J(X))$ be
  arbitrary. Since the group $B$ is Boolean, $\sigma$ may be written
  as $\sigma=\gamma_1\gamma_2\ldots \gamma_s$ where the elements
  $\gamma_1, \gamma_2,\ldots, \gamma_s$ are different generators of
  $B$, and this decomposition is unique up to reordering of the
  factors. Since $\sigma \in \xi(\Pi^x_I(X))$, and because the voltage
  of a semi-edge is always a generator, each $\gamma_i$ is also in
  $\xi(\Pi^x_I(X))$, and since $\sigma \in \xi(\Pi^x_J(X))$ each
  $\gamma_i$ is also in $\xi(\Pi^x_J(X))$. But, this implies that each $\gamma_i$ is in
  $\xi(\Pi^x_{I\cap J}(X))$, implying that
  $\sigma \in \xi(\Pi^x_{I\cap J}(X))$. Therefore,
  $\xi(\Pi^x_I) \cap \xi(\Pi^x_J) = \xi(\Pi^x_{I \cap J})$. \qed
\end{dem}

Proposition~\ref{p:Caterpillars1} is still true for infinite
caterpillars. Even if our algorithmic way of assigning voltages may
not be doable, the voltage assignment is still well defined as a
quotient of the Boolean group with one generator assigned to each
semi-edge.

Note that Proposition~\ref{p:Caterpillars1} says that every caterpillar $X$ is a symmetry type graph of a polytope $\p$ with respect to some group $B\leq \Aut(\p)$.
However need not happen that $X$ is the symmetry type graph of $\p$ with respect to the full automorphism group.
That is, $X^\xi$ might have ``extra'' symmetry.
Thus, we would want to investigate what could be the
symmetry type of $X^\xi$ with respect to its full automorphism group.

Let $X$ be a finite caterpillar and let $\xi:\Pi(X)\to B$ be the
voltage assignment constructed in the proof of
Proposition~\ref{p:Caterpillars1}. If $X$ is symmetric, its
non-trivial symmetry induces an automorphism of $B$ which is just a
reordering of the generators, moreover, if the generators are given
the natural order, it induces the reverse order on the
generators. Using~\cite[Theorem 7.1]{Voltajes}, we get that this symmetry induces a symmetry of the derived maniplex, so in this case the original caterpillar is
not the symmetry type of the derived polytope. But we will see in
Theorem~\ref{t:Caterpillars} that if this is not the case we can
almost be certain that the caterpillar is in fact the symmetry type
graph of the derived polytope by its full automorphism group. In this
case, by ``almost'' we mean that if this is not the case, the
caterpillar must have a very specific structure.

\begin{teo}\label{t:Caterpillars}
  Let $X$ be a finite caterpillar of length $k$ and rank $n$. Let
  $S=c_1c_2\ldots c_k$ be the underlying word of $X$. Then at
  least one of the following statements is true:
  \begin{enumerate}
  \item $X$ is symmetric.
  \item $X$ is the STG of a polytope with a Boolean automorphism
    group.
  \item $c_1$ is in $\{1,n-2\}$ and there exist $r\in
    \{1,2,\ldots,k-1\}$, and
    $a_1,a_2,\ldots,a_t \in \{0,1,\ldots,n-1\}$ where
    $t=(k+1)/(2r+2)$, such that
    $S=wa_1w^{-1}bwa_2w^{-1}b\ldots bwa_tw^{-1}$ where $w=c_1c_2\ldots
    c_r$ and 
    $b=0$ if $c_1=1$ and $b=n-1$ if $c_1=n-2$.
  \item There exist $r\in \{1,2,\ldots,k-1\}$ and
    $a,b \in \{0,n-1\}$ such that
    \[
      S=waw^{-1}bwaw^{-1}b\ldots bwaw^{-1}bw,
    \]
    where $w=c_1c_2\ldots c_r$. Also
    $(c_1,b), (c_r,a) \in \{(1,0),(n-2,n-1)\}$.
  \end{enumerate}
\end{teo}
\begin{dem}
  Suppose that $X$ is not symmetric and that it is not the STG of a
  polytope with a Boolean automorphism group. 

  Consider the voltage
  assignment $\xi:\Pi(X)\to B$ previously discussed. We say that
  two vertices $x$ and $y$ of $X$ are equivalent ($x\sim y$) if there
  exist (or equivalently, for all) $\sigma,\tau \in B$ such that
  the flags $(x,\sigma)$ and $(y,\tau)$ of $X^\xi$ are in the same
  orbit under the action of the automorphism group of $X^\xi$. Then
  $\sim$ is an equivalence relation preserved by $i$-adjacency,
  that is $x\sim y \Rightarrow x^i\sim y^i$. Moreover,
  the natural function $h:X\to X/\sim$ is a premaniplex
  homomorphism, so by Proposition~\ref{p:CubsCat} there
  exists some $r$ such that the $S$ can be written
  as $wa_1w^{-1}b_1wa_2w^{-1}b_2\ldots $ ending after an occurrence of
  either $w$ or $w^{-1}$, where $w=c_1c_2\ldots c_r$, $a_i\in
  \{c_r+1,c_r-1\}$ and $b_i\in \{c_1+1,c_1-1\}$ (see
  Figure~\ref{f:FoldedCat}).


  Note that since $X$
  is not symmetric nor the STG of $X^\xi$ with respect to $B$, at
  least $b_1$ exists. Let $j\in\{0,1\ldots,k-1\}$ be a number such
  that the segment $[x_0,x_{j+1}]$ has the underlying word
  $wa_1w^{-1}b_1wa_2w^-1b_2\ldots w^{-1}b_i$ for some $i$. We know in
  particular that $b_i$ differs from $c_1$ in exactly 1. We want to
  prove that $(c_1,b_i) \in \{(1,0),(n-2,n-1)\}$. We will assume this
  is not the case and arrive to a contradiction.

  Let $q$ be the {\em other} color that differs from
  $b_i$ in exactly 1 (that is $q = 2b_i - c_1$). Since $(c_1,b_i) \neq
  (1,0),(n-2,n-1)$ we know that $q\in \{0,1,\ldots,n-1\}$, and thus is
  the color of some edges of $X$. So there
  are semi-edges $e,e'$ incident
  to $x_0$ of colors $q$ and $b_i$ respectively. Let
  $\alpha:=\xi(e)$ and $\beta:=\xi(e')$. The voltage of the closed
  path $ee'ee'$ is $(\beta\alpha)^2 = 1$ because $B$ is Boolean. This
  means that its lift, (the path of length 4 that starts at $(x,1)$ in
  $X^\xi$ and  alternates colors between $q$ and $b_i$) must be closed.

  %


  By Theorem~\ref{p:CubsCat}, we know that $c_{j+1} = c_1 \neq q$, so
  we know that the darts $(x_{j+1},q)$ and $(x_j,q)$ are
  semi-edges. Let $\kappa := \xi(x_j,q)$ and
  $\lambda = \xi(x_{j+1},q)$. The path of length 4 that alternates
  colors between $r$ and $b_i$ and starts at $x_j$ is closed, and its
  voltage is $\lambda\kappa$. Note that since $|q-b_i|=1$ the
  construction of $\xi$ tells us that
  $\xi(x_j,q) \neq \xi(x_{j+1},q)$, that is $\lambda \neq \kappa$,
  which implies $\lambda\kappa \neq 1$. This means that the path of
  length 4 in $X^\xi$ starting at $(x_j,1)$ and alternating colors
  between $r$ and $b_i$ is not closed (it ends at
  $(x_j,\lambda\kappa)$).

    \begin{figure}
    \begin{center}
      \includegraphics[width=7cm]{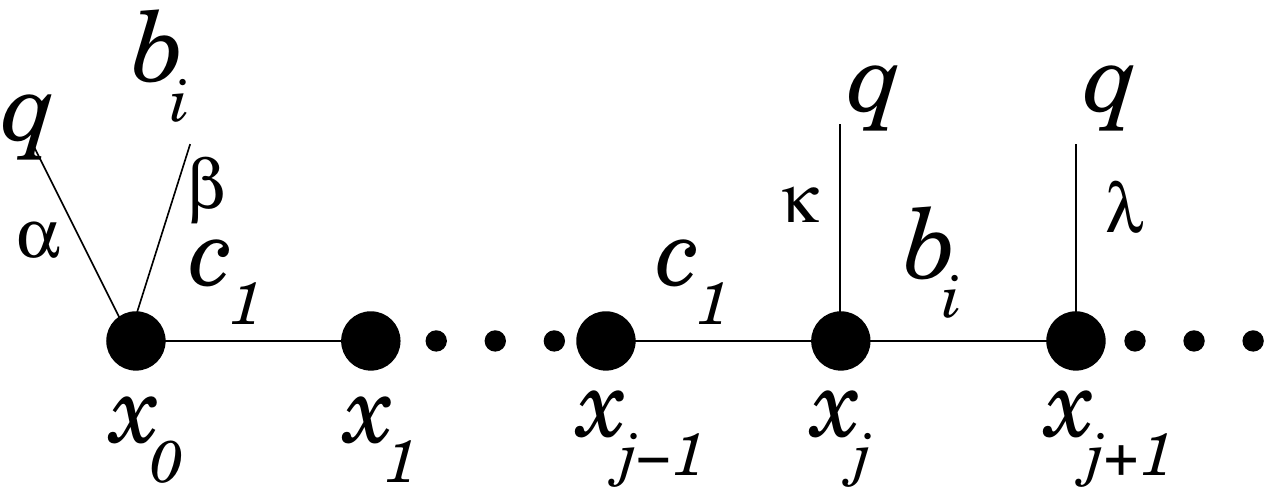}
      \caption{\label{f:CatCaso2} 
        The voltage of the path that alternates colors between $r$ and
        $b_i$ starting at $x_j$ is $\lambda\kappa\neq 1$.}
    \end{center}
  \end{figure}

  We see that the path of length 4 in $X^\xi$ starting at $(x_j,1)$
  and alternating colors between $r$ and $c_{j+1}$ is not closed, but
  the one starting at $(x_0,1)$ is. This contradicts the fact that
  $(x_j,1)$ and $(x_0,1)$ are on the same orbit. The contradiction
  comes from the fact that there are edges of color
  $q = 2b_i -c_1 \in \{0,1,\ldots,n-1\}$, so to avoid this we must
  have that $(c_1,b_i) \in \{(1,0),(n-2,n-1)\}$. Since $c_1$ is fixed,
  every $b_i$ must be the same.


  If the underlying word of $X$ ends after an occurrence of $w$ we
  may look at $X$ in the other direction. Then the previous result
  tells us that every $a_i$ is equal to some $a$ and that $(c_r,a)\in
  \{(1,0), (n-2,n-1)\}$. \qed
\end{dem}

 \begin{obs}
   If the third or fourth condition is the one holding in
   Theorem~\ref{t:Caterpillars}, the actual STG of $X^\xi$ might be the
   finite caterpillar with underlying word $w$.
 \end{obs}

By doing exactly the same proof, we obtain the following analogous
result for infinite caterpillars with one end-point (this would
  look like a ray or half-straight line):

\begin{teo}\label{t:Caterpillars1end}
  Let $X$ be an infinite caterpillar with one end-point. Let $S$ be
  the sequence of colors of the underlying path of $X$ starting at its
  end-point. Then one of the following statements is true:
  \begin{enumerate}
  \item $X$ is the STG of a polytope with a Boolean automorphism
    group.
  \item There exist some number $r$ and colors
    $b,a_1,a_2,\ldots,a_t \in \{0,1,\ldots,n-1\}$ such that
    $S=wa_1w^{-1}bwa_2w^{-1}bwa_3w^{-1}\ldots$ where $w=c_1c_2\ldots
    c_r$ and $(c_1,b) \in \{(1,0),(n-2,n-1)\}$.
  \end{enumerate}
\end{teo}

 \begin{obs}
   If the second condition is the one holding, the actual STG of
   $X^\xi$ might be the finite caterpillar with underlying word $w$.
 \end{obs}

If $k\geq 3$, it is easy to construct a caterpillar of length $k-1$ (that is, with $k$ vertices and $k-1$ links) that does not satisfy properties $(1)$, $(3)$ and $(4)$ (simply let $c_1=0$ and then avoid the color $0$ for every other $c_i$).
Hence, we get the following corollary:

\begin{coro}\label{coro:k_orbit_exists}
    For every $n,k\geq 3$, there is an abstract $n$-polytope with Boolean automorphism group that has $k$ flag orbits.
\end{coro}



\section*{Acknowledgements}
The authors would like to thank Gabe Cunningham, Antonio Montero and Daniel Pellicer for helpful comments on early versions of this paper. We also thank the financial support of CONACyT grant A1-S-21678.

\bibliographystyle{plain}
\bibliography{Bibliografia}

\end{document}